\documentclass{article}
\usepackage{arxiv}

\usepackage[utf8]{inputenc} 
\usepackage[T1]{fontenc}    
\usepackage[hidelinks]{hyperref}     
\usepackage{url}            
\usepackage{booktabs}       
\usepackage{amsfonts}       
\usepackage{nicefrac}       
\usepackage{microtype}      
\usepackage{lipsum}
\usepackage{bm}
\usepackage{graphicx}
\usepackage{subcaption}
\usepackage{amsmath}
\graphicspath{ {./images/} }

\title{Transient thermal analysis of a bi-layered composites with the dual-reciprocity inclusion-based boundary element method}

\author{
 Chunlin Wu \\
  Shanghai Institute of Applied Mathematics and Mechanics,\\
  Shanghai University\\
  Shanghai, 200044, China \\
  \texttt{chunlinwu@shu.edu.cn} 
    \And
  Liangliang Zhang \\
  Department of Applied Mechanics, \\China Agricultural University, \\Beijing, 100083, China\\
  \texttt{llzhang@cau.edu.cn}
   \And
 Tengxiang Wang \\
 Division of Engineering,\\ Business and Computing, \\Penn State Berks, \\Reading, PA, 19610, United States\\
  \texttt{tbw5524@psu.edu} 
  \And
 Huiming Yin \\
  Department of Civil Engineering  and Engineering Mechanics,\\
  Columbia University\\
  New York, NY, 10027, United States \\
  \texttt{yin@civil.columbia.edu} \\
}

\begin{document}
\maketitle
\begin{abstract}
This paper proposes a single-domain dual-reciprocity inclusion-based boundary element method (DR-iBEM) for a three-dimensional fully bonded bi-layered composite embedded with ellipsoidal inhomogeneities under transient/harmonic thermal loads. The heat equation is interpreted as a static one containing time- and frequency-dependent nonhomogeneous source terms, which is similar to eigen-fields but is transformed into a boundary integral by the dual-reciprocity method. Using the steady-state bimaterial Green's function, boundary integral equations are proposed to take into account continuity conditions of temperature and heat flux, which avoids setting up any continuity equations at the bimaterial interface. Eigen-temperature-gradients and eigen-heat-source are introduced to simulate the material mismatch in thermal conductivity and heat capacity, respectively. The DR-iBEM algorithm is particularly suitable for investigating the transient and harmonic thermal behaviors of bi-layered composites and is verified by the finite element method (FEM). Numerical comparison with the FEM  demonstrates its robustness and accuracy. The method has been applied to a functionally graded material as a bimaterial with graded particle distributions, where particle size and gradation effects are evaluated.
\end{abstract}

\keywords{Inhomogeneity problem \and Time-harmonic/transient thermal analysis \and  Bimaterial Green's function \and  Inclusion-based boundary element method \and  Functionally graded material}

\section{Introduction}
The composite material with two or more dissimilar material phases plays an important role in engineering applications. For example, the bi-layered photovoltaic panels in civil engineering improve overall energy harvest efficiency \cite{zadshir2023}, and the functionally graded materials in combustion engines significantly reduce thermal stresses between two dissimilar material (ceramic and metal) interfaces \cite{saleh2020}. Analytical and numerical methods have been applied to understand and optimize the heat transfer process. Green's functions serve as particular solutions of the governing equation, which relates excitations and field variables under specified boundary conditions. Green's functions for bi-materials are very useful as they can be reduced to a single infinite domain and semi-infinite domain by changing the material constants \cite{Wang_2022_JAM,Yin_iBEM}. 

In the literature, different Green's functions have been proposed for the different material domains. For instance, Mindlin's solution \cite{Mindlin1936}, Walpole's solution \cite{Walpole1996}, and Yue's solution \cite{Yue1995} provided the displacement field caused by a concentrated force located in the semi-infinite, bimaterial, and multi-layered space, respectively. Their pioneering works with the method of images paved the way for later prosperity of Green's function for much more complicated cases, such as transversely isotropic semi-infinite in geomechanics \cite{Pan1979}, and certain distribution of surface loads \cite{Xiao_2011}. Several works in the literature were devoted to Green's functions for steady-state heat conduction. Berger et al. \cite{Berger2000} proposed a steady-state Green's function for heat conduction in anisotropic bimaterial, and the authors employed the virtual force method to express homogeneous solutions. Ang et al. \cite{Ang2004} considered the two-dimensional Green's function, which considers an imperfect bimaterial interface with discontinuous interfacial heat flux. Subsequently, Wang et al. \cite{Wang2007} proposed a Green's function for a graded half-space, and Li et al. \cite{Li2020} derived solutions for a bimaterial problem with fluid and isotropic solids. However, it is challenging to extend the method of Green's function from the steady-state case to the transient case. Zhou et al. \cite{Zhou2020} derived the three-dimensional Green's function for anisotropic bimaterial, which was expressed in the Fourier/Laplace space instead. 
 
Since the direct application of Green's function has strict limitations on either boundary or geometric conditions, numerical methods have become popular alternative tools for thermal analysis, such as the boundary element method (BEM) and finite element method (FEM), and their numerous extensions. Particularly, BEM relies on Green's functions with their reciprocal properties and requires boundary discretization for boundary integrals with Green's function. Thanks to the rapid development of Green's function for several boundary conditions, such as semi-infinite, bimaterial \cite{Wang_2022_JAM}, the seamless switch from the infinite Green's functions to bimaterial Green's function exhibits unique advantages. Recently, Wu et al. \cite{Wu2023_EABE, Wu2024_IJES} used the bimaterial elastic/thermoelastic Green's function as the kernel functions, and developed a single-domain boundary element method for bi-layered composites. Since the bimaterial Green's function mathematically involves the continuity conditions across the bimaterial interface, two primary merits are obtained: (i) no interface discretization is required; therefore, the linear system contains fewer degrees of freedom; (ii) a more accurate representation of interior points, as the ill-conditioned integral of kernel functions in the neighborhood of the bimaterial interface has been avoided. However, bimaterial Green's functions for transient problems \cite{Zhou2020} were often given in implicit or numerical forms or transformed spaces, which limit the accuracy and efficiency in application.

Pioneering works in the dual-reciprocity methods shed light on the replacement of these Green's functions in the transient state by the counter-parts in the steady-state \cite{Nardini1983}. Partridge et al. \cite{Partridge1991, Partridge2000} proposed to regard the original transient equations as a nonhomogeneous quasi-static equation \cite{Zhou2012}. In analogy with elastodynamics, the heat generation rate term in the transient heat equation can be considered as an inhomogeneous heat source. Although the evaluation of heat source requires domain integral of Green's function, it can be converted into boundary integrals using the reciprocal property of Green's function twice \cite{Kamiya1993}. Ang et al. \cite{Ang2006} proposed a two-dimensional dual-reciprocity boundary element method (DRBEM) with the implementation of bimaterial steady-state Green's function \cite{Ang2004}, which exhibited good agreement with the analytical solutions. Following their concepts, this paper applies the three-dimensional DRBEM with bimaterial Green's function provided in \cite{Wang_2022_JAM} to the inhomogeneity problems. Generally, when the bi-layered composites contain with many inhomogeneities, the DRBEM still requires interface discretization among inhomogeneity surfaces, following the continuity of the multi-domain BEM \cite{Beer2008}. However, using Eshelby's equivalent inclusion method (EIM), we can avoid the multi-domain boundary integral. 

Recently, we integrated the DRBEM and EIM into a new method, namely the dual-reciprocity (DR-iBEM). While the former handles various boundary responses, the latter allows the simulation of inhomogeneities without surface mesh. The EIM for transient/harmonic heat transfer is a subsequent extension of the steady-state one \cite{Eshelby_1957, Eshelby_1959} following the concept of the polarization method \cite{willis1997dynamics}. Specifically, the equivalent inclusions exhibit the same material properties as the matrix, and the thermal conductivity and heat capacity mismatch are simulated by eigen-temperature-gradient (ETG) \cite{hatta1986} and eigen-heat-source (EHS) \cite{Wu2024-prsa}, respectively. Since DRBEM utilizes steady-state Green's functions, disturbances caused by two eigen-fields within inhomogeneities are evaluated through closed-form domain integrals of steady-state Green's functions, which are known as steady-state Eshelby's tensors \cite{Wang_2022_JAM, Wu2024-prsa}. Note that the DR-iBEM framework avoids the need to derive transient/dynamic Eshelby's tensors defined by transient Green's functions. Therefore, all inhomogeneities can be handled by the well-established domain integrals in the literature, such as elliptical/ellipsoidal \cite{Moschovidis1975}, and polygonal/polyhedral \cite{Wu2021_jam_polygonal, Wu2021_polyhedral} inclusions. 

This paper develops the dual-reciprocity inclusion-based boundary element method (DR-iBEM) to bi-layered composites containing multiple ellipsoidal inhomogeneities, which is a single-domain numerical method for transient/harmonic thermal analysis with the interfacial continuity assured by the Green's function. Section 2 introduces the boundary value problem of a bi-layered composite containing inhomogeneities with continuity equations and eigen-fields. Section 3 proposes the three-dimensional (3D) DRBEM with bimaterial Green's function, and disturbances caused by eigen-fields are determined by equivalent conditions. Section 4 verifies the DR-iBEM by comparing the local thermal fields obtained from FEM. Section 5 conducts transient and harmonic thermal analysis of a functionally graded material (FGM). Finally, some conclusive remarks are provided. 

\section{Problem statement}
\begin{figure}
    \centering
    \includegraphics[width=0.4\linewidth]{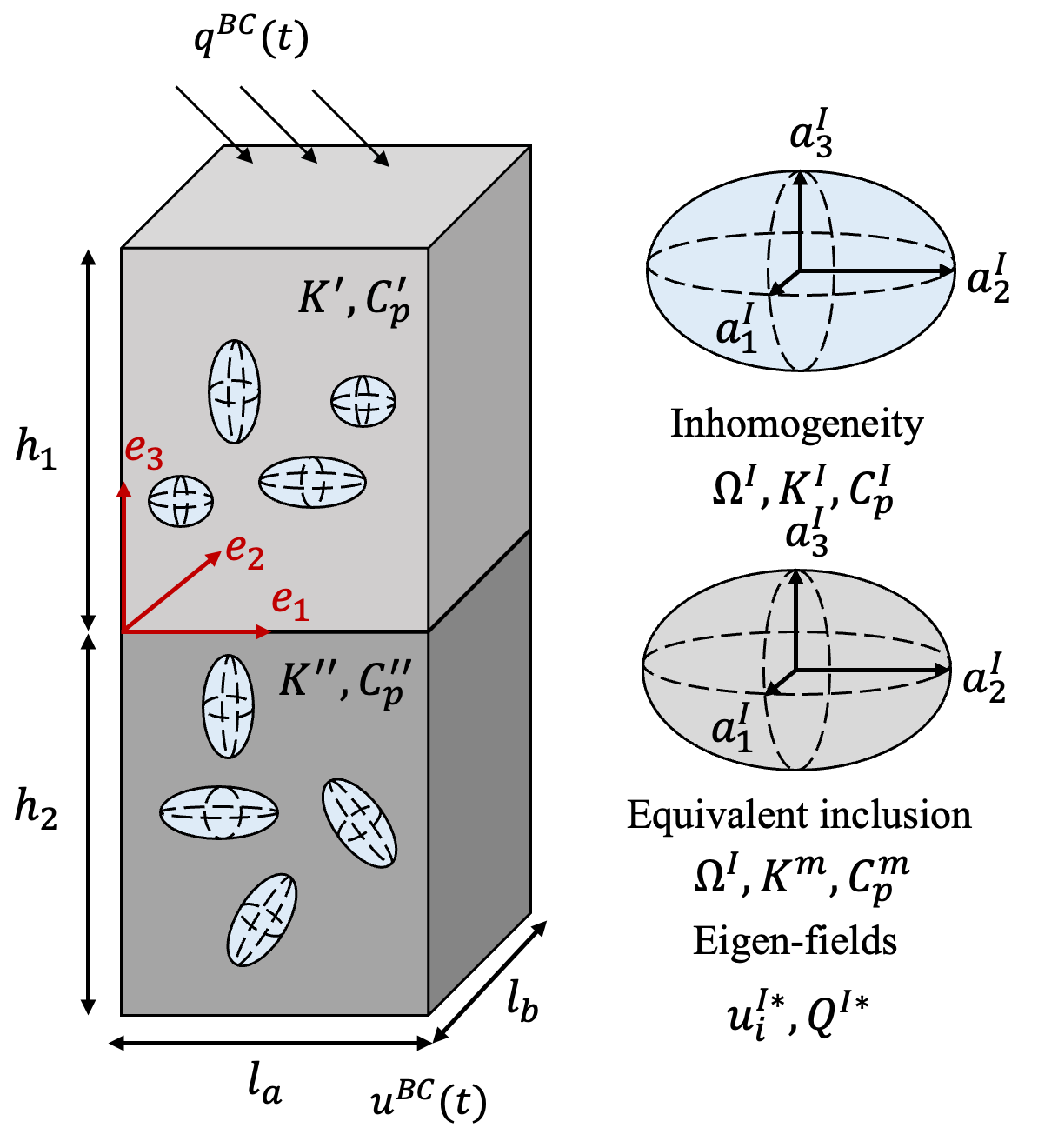}
    \caption{Schematic illustration of a bi-layered composites system composed of two dissimilar jointed blocks with same length $l_a$ and width $l_b$, but different heights at $h_1$ and $h_2$, respectively. Multiple ($N^I$) ellipsoidal subdomains $\Omega^I$ are embedded in the bi-layered composites, which exhibit thermal conductivity $K^I$ and specific heat $C_p^I$, and they can be handled as equivalent inclusions with continuously distributed ETG and EHS.}
    \label{fig:geometryc}
\end{figure}

Fig. \ref{fig:geometryc} schematically illustrates a bi-layered composite sample containing multiple ($N^I$) ellipsoidal inhomogeneities ($\Omega^I$), which is subjected to prescribed temperature ($u^{BC}$) and heat flux ($q^{BC}$) on the boundary. The dimension of the bi-layered composite is length $l_a$, width $l_b$, and the upper and lower blocks exhibit the height $h_1, h_2$, respectively. Each ellipsoidal domain is characterized by its center ($\textbf{x}^{IC}$) and three radii $\textbf{a}^I = (a_1^I, a_2^I, a_3^I)$. Generally, two matrix phases and inhomogeneities exhibit different thermal conductivity and specific heat constants as $K$ and $C_p$, respectively. Specifically, the thermal properties belonging to the upper, lower matrix, and inhomogeneities are denoted by superscripts $(.)',(.)'', (.)^I$, respectively. 

This paper assumes two matrix phases and embedded inhomogeneities are perfectly bonded, and the temperature and flux satisfy continuity equations along their interfaces, including the bimaterial interface ($x_3 \equiv 0$). Therefore, a well-posed boundary value problem is formed with the continuity conditions and boundary conditions for the bi-layered composite. To avoid any interior discretization or additional unknowns for continuity conditions, this paper adopts two steps: (i) each inhomogeneity is treated as an equivalent inclusion, exhibiting the same thermal properties as the surrounding matrix with eigen-fields to simulate the material mismatch; (ii) disturbances by eigen-fields within equivalent inclusions and continuity conditions of bimaterial interface are taken into account by Green's function through domain integral on inclusions and boundary integral on the outer boundary. Note that Green's function handles the inhomogeneity surface and bi-material interface continuity automatically.

Readers may question why the same iBEM algorithm is not extended using the bimaterial transient Green's function, particularly since it can be reduced from Zhou's solution \cite{Zhou2020} from the general anisotropic case. However, it is essential to note that Zhou and Han \cite{Zhou2020} derived the bimaterial transient Green's function in the transformed domain, and there exist no closed-form formulae in the temporal-spatial domain, which require the two-dimensional inverse Fourier and inverse Laplace transforms. These inverse numerical transforms inevitably increase computational burdens, and subsequent implementations with boundary integrals become more time-consuming and involve numerical accuracy issues, undermining Green's function's primary advantage. Moreover, similar problems apply to transient Eshelby's tensors, which involve temporal-spatial domain integrals of Green's function. To avoid such issues in numerical transform, this paper alternatively adopts the dual-reciprocity method, which allows transient thermal analysis to be performed with the steady-state bimaterial Green's function and closed-form Eshelby's tensors.

\section{Formulation}

\subsection{Bimaterial dual-reciprocity boundary integral equations}

In the previous section, the DR-iBEM algorithm handles all inhomogeneities as equivalent inclusions, which converts the original bi-layered heterogeneous system into an equivalent bi-layered homogeneous system with source fields distributed over the subdomains. In the absence of a prescribed heat source, the heat equation for the bi-layered homogeneous system can be written as Eq. (\ref{eq:heat_eqn}):

\begin{equation}
    K(\textbf{x}) u_{,ii}(\textbf{x}, t) = C_p(\textbf{x}) \frac{\partial u(\textbf{x}, t)}{\partial t}
    \label{eq:heat_eqn}
\end{equation}
where $C_p(\textbf{x}) \frac{\partial u(\textbf{x}, t)}{\partial t}$ can be treated the same as an inhomogeneous heat source. It should be noted that the thermal conductivity and heat capacity are piecewise constant functions depending on the position of $x_3$. Specifically, (i) $x_3 \geq 0$, $C_p = C_p'$, $K = K'$; otherwise (ii) $x_3 < 0$, $C_p = C_p''$, $K = K''$. In such a case, when the heat equation is interpreted as a static one with a time-dependent inhomogeneous heat source, the bimaterial Green's function has been derived in our recent work \cite{Wang_2022_JAM}:
\begin{equation}
    G(\textbf{x}, \textbf{x}') = \begin{cases} \frac{1}{4 \pi K^s} \phi + \frac{1}{4\pi K^s} \frac{K^s - \overline{K}^s}{K^s + \overline{K}^s} \overline{\phi} & x_3 x_3' \geq 0 \\ \frac{1}{2\pi (K^s + \overline{K}^s)} \phi & x_3 x_3' < 0 \end{cases} 
    \label{eq:Gfunc}
\end{equation}
where $K^s$ denotes the thermal conductivity of the material at the source point $\textbf{x}'$ and $\overline{K}^s$ the thermal conductivity at the image point $\overline{\textbf{x}}'$, respectively, which depend on the relative position of field and source points as follows: (i) $x_3' < 0$, $K^s ='$, $\overline{K}^s = ''$; and (ii) $x_3' \geq 0$, $K^s = ''$, $\overline{K}^s = '$; $\quad \phi = \frac{1}{|\textbf{x} - \textbf{x}'|} \quad  \text{and} \quad \overline{\phi} = \frac{1}{|\textbf{x} - \overline{\textbf{x}}'|}$, in which  $\overline{\textbf{x}}'$ is the imaged source point as $\overline{x}_i' = M_I x_i'$ with $\textbf{M} = \{1, 1, -1\}$, where the dummy index summation does not apply to the capital index $I$ \cite{Mura1987}, while $I$ changes with $i$. Although the bimaterial Green's function contains both $\phi$ and $\overline{\phi}$, it satisfies the governing equation:
\begin{equation}
    -K(\textbf{x}) G_{,ii}(\textbf{x}, \textbf{x}') = \delta(\textbf{x}-\textbf{x}')
\end{equation}
where $\delta(\textbf{x}-\textbf{x}')$ is the Dirac Delta function. Although boundary integral equations for steady-state heat conduction in a bi-layered system have been proposed in \cite{Wu2024_IJES}, it is important to formulate the boundary integral equations for transient heat transfer and demonstrate their primary differences. Using the filtering effect of the Dirac Delta function, the thermal fields at any interior points within the domain $\mathcal{D}$ can be expressed as:
\begin{equation}
\begin{split}
    &c(\textbf{x}) u(\textbf{x}, t)  = \int_{\mathcal{D}} \delta(\textbf{x}, \textbf{x}') u(\textbf{x}', t) \thinspace d\textbf{x}' = \int_{\mathcal{D}} -K(\textbf{x}') G_{,ii}(\textbf{x}, \textbf{x}') u(\textbf{x}', t) \thinspace d\textbf{x}' \\ 
    & = \int_{\mathcal{D}} \frac{\partial}{\partial x_i'} \left[ -K(\textbf{x}') G_{,i'}(\textbf{x}, \textbf{x}') u(\textbf{x}', t) \right] \thinspace d\textbf{x}' - \int_{\mathcal{D}} \frac{\partial }{\partial x_i'} \left[ -K(\textbf{x}') u(\textbf{x}', t) \right] G_{,i'}(\textbf{x}, \textbf{x}') \thinspace d\textbf{x}' \\
    & = -\int_{\partial \mathcal{D}} n_i(\textbf{x}') K(\textbf{x}') G_{,i'}(\textbf{x}, \textbf{x}') \thinspace d\textbf{x}' + \int_{\mathcal{D}^+} K' u_{,i'}(\textbf{x}, t) G_{,i'}(\textbf{x}, \textbf{x}') \thinspace d\textbf{x}' + \int_{\mathcal{D}^-} K'' u_{,i'}(\textbf{x}, t) G_{,i'}(\textbf{x}, \textbf{x}') \thinspace d\textbf{x}' \\
    & = -\int_{\partial \mathcal{D}} n_i(\textbf{x}') K(\textbf{x}') G_{,i'}(\textbf{x}, \textbf{x}') \thinspace d\textbf{x}' + \int_{\mathcal{D}^+} K' \frac{\partial }{\partial x_i'} \left[ u_{,i'}(\textbf{x}', t) G(\textbf{x}, \textbf{x}') \right] +  \int_{\mathcal{D}^-} K'' \frac{\partial }{\partial x_i'} \left[ u_{,i'}(\textbf{x}', t) G(\textbf{x}, \textbf{x}') \right] \thinspace d\textbf{x}' \\ 
    & \quad \thinspace - \int_{\mathcal{D}^+} K' u_{,i'i'}(\textbf{x}', t) G(\textbf{x}, \textbf{x}') \thinspace d\textbf{x}'- \int_{\mathcal{D}^-} K'' u_{,i'i'}(\textbf{x}', t) G(\textbf{x}, \textbf{x}') \thinspace d\textbf{x}' \\
    & = -\int_{\partial \mathcal{D}} n_i(\textbf{x}') K(\textbf{x}') G_{,i'}(\textbf{x}, \textbf{x}') \thinspace d\textbf{x}' + \int_{\partial \mathcal{D}^+} K' n_i(\textbf{x}') u_{,i'}(\textbf{x}, \textbf{x}') G(\textbf{x}, \textbf{x}') \thinspace d\textbf{x}' \\ 
    & \quad \thinspace + \int_{\partial \mathcal{D}^-} K'' n_i(\textbf{x}') u_{,i'}(\textbf{x}, \textbf{x}') G(\textbf{x}, \textbf{x}') \thinspace d\textbf{x}' - \int_{\mathcal{D}^+} C_p' G(\textbf{x}, \textbf{x}') \frac{\partial u(\textbf{x}', t)}{\partial t} \thinspace d\textbf{x}'- \int_{\mathcal{D}^-} C_p'' G(\textbf{x}, \textbf{x}') \frac{\partial u(\textbf{x}', t)}{\partial t} \thinspace d\textbf{x}'
\end{split}
    \label{eq:BIE_source}
\end{equation}
where $c(\textbf{x})$ is the free term, which is $\frac{1}{2}$ and $1$ for smooth boundary and interior points, respectively \cite{Beer2008}; $\textbf{n}$ represent the unit outward normal vectors on the boundary point $\textbf{x}'$; $\partial \mathcal{D}$ refers to the boundary of the domain $\mathcal{D}$; $K(\textbf{x}') u_{,i'}(\textbf{x}', t) = q^{BC}(\textbf{x}', t)$ is the boundary ''flux'', and; $n_i$ is the unit outward normal of the boundary point. Note that since the bimaterial Green's function naturally involves continuity conditions of temperature and flux, two boundary integral terms can be united as one exterior boundary integral: 
\begin{equation}
    \int_{\partial \mathcal{D}^+} K' n_i(\textbf{x}') u_{,i'}(\textbf{x}, \textbf{x}') G(\textbf{x}, \textbf{x}') \thinspace d\textbf{x}' + \int_{\partial \mathcal{D}^-} K'' n_i(\textbf{x}') u_{,i'}(\textbf{x}, \textbf{x}') G(\textbf{x}, \textbf{x}') \thinspace d\textbf{x}' = \int_{\partial \mathcal{D}} K(\textbf{x}) n_i(\textbf{x}') u_{,i'}(\textbf{x}') G(\textbf{x}, \textbf{x}') \thinspace d\textbf{x}'
\end{equation}
where the boundary integrals along the interface of the bimaterial are canceled due to opposite normal vectors and continuous Green's function across the interface. Hence, Eq. (\ref{eq:BIE_source}) can be subsequently simplified:
\begin{equation}
\begin{split}
    c(\textbf{x}) u(\textbf{x}, t) & + \int_{\partial \mathcal{D}} \frac{\partial G(\textbf{x}, \textbf{x}')}{\partial x_i'} K(\textbf{x}') n_i(\textbf{x}') u^{BC}(\textbf{x}', t) \thinspace d\textbf{x}' \\ &= \int_{\partial \mathcal{D}} G(\textbf{x}, \textbf{x}') q^{BC}(\textbf{x}', t) \thinspace d\textbf{x}' - \int_{\mathcal{D}^+} C_p' G(\textbf{x}, \textbf{x}') \frac{\partial u(\textbf{x}', t)}{\partial t} \thinspace d\textbf{x}'- \int_{\mathcal{D}^-} C_p'' G(\textbf{x}, \textbf{x}') \frac{\partial u(\textbf{x}', t)}{\partial t} \thinspace d\textbf{x}'
\end{split}
    \label{eq:BIE_con}
\end{equation}
Note that the last two domain integrals in Eq. (\ref{eq:BIE_con}) cannot be combined as a single one, unless the upper and lower materials exhibit exactly the same heat capacity. In addition, the bimaterial Green's function only ensures the continuity of temperature and heat flux across the bimaterial interface, and different heat capacities multiplied by temperature are not continuous across the bimaterial interface. An inappropriate combination of these two domain integrals will lead to discontinuous temperature and flux distribution across the bimaterial interface. 

The heat generation/storage rate $C_p(\textbf{x}) \frac{\partial u(\textbf{x}, t)}{\partial t}$ involves the partial derivative of temperature with respect to time. The time derivatives can be approximated by the second-order Euler scheme, while the first-order scheme is applied for the first time station: 
\begin{equation}
    \frac{\partial u(\textbf{x}, t_n)}{\partial t} =
    \begin{cases}
    \frac{3 u(\textbf{x}, t_n) - 4 u(\textbf{x}, t_{n-1}) + u(\textbf{x}, t_{n-2})}{2 \Delta t} & n > 1 \\
    \frac{u(\textbf{x}, t_1) - u(\textbf{x}, t_0)}{\Delta t} & n = 1
    \end{cases}
    \label{eq:backward_temp}
\end{equation}
where $t_0 = 0$ is the initial time station; $t_n = t_0 + n \Delta t$ is the $n^{th}$ time station (n = 0, 1, 2, ...), and $\Delta t$ is the time interval. For transient heat transfer, this paper only presents formulae for $n \geq 2$, because the formulae can be straightforwardly reduced to the case for the first time station. 

When the system is subjected to a sinusoidal thermal load,  the response can reach the time-harmonic state so that the spatial and temporal variables can be separated, which is potentially applicable to the seasonal temperature variation in geothermal engineering. Specifically, the temperature field can be separated as $u(\textbf{x}, t) = \tilde{u}(\textbf{x}) \exp[-i \omega t]$, where $\omega$ is the circular frequency of thermal loads. Therefore, the time partial derivative can be derived as follows: 
\begin{equation}
    \frac{\partial u(\textbf{x}, t)}{\partial t} = -i \omega u(\textbf{x}, t)
    \label{eq:dif_harmonic}
\end{equation}

The dual-reciprocity method (DRM) has been presented in several papers \cite{Partridge1991, Partridge2000, Ang2006}, which applied the reciprocal properties of Green's function to convert domain integrals into surface integrals. Specifically, the temperature field is approximated through the nodal values and radial basis functions (RBF). Since this paper aims to convey the dual-reciprocity iBEM algorithm, the approximation scheme is inherited from \cite{Ang2006}, and only necessary derivations for the bimaterial boundary integral equations are kept. At the time station $t_n$, the temperature distribution of any interior point can be expressed in terms of the superposition of $NN$ boundary nodes and $NS$ interior nodes: 
\begin{equation}
\begin{aligned}
    u(\textbf{x}, t_{n}) = \sum_{m=1}^{NN + NS} \chi(\textbf{x}, \textbf{x}^m) \thinspace \alpha(\textbf{x}^m, t_{n}) 
\end{aligned}
    \label{eq:interpolation_temp} 
\end{equation}
where $\alpha(\textbf{x}^m, t_n)$ refers to the virtual source density at the time station $t_n$, which is created to satisfy the temperature distribution using the present interpolation with RBF for $NN + NS$ collocation points; RBF $ \chi(\textbf{x}, \textbf{x}^m) = 1 + |\textbf{x} - \textbf{x}^m| + |\textbf{x} - \textbf{x}^m|^{2}$, which defines $\nabla^2 \Gamma(\textbf{x}, \textbf{x}^m) =\chi(\textbf{x}, \textbf{x}^m)$. Therefore $\Gamma(\textbf{x}, \textbf{x}^m)$ 
can be explicitly expressed as $\Gamma(\textbf{x}, \textbf{x}^m) = \frac{|\textbf{x} - \textbf{x}^m|^2}{6} + \frac{|\textbf{x} - \textbf{x}^m|^3}{12} + \frac{|\textbf{x} - \textbf{x}^m|^4}{20}$. Note that the selection of RBF is based on the recommendation in \cite{Partridge2000}. However, different forms of RBF can be used \cite{Ang2006}, which leads to different $\Gamma(\textbf{x}, \textbf{x}^m)$. For the steady-state response under a sinusoidal load, the spatial part can be separated from the temporal part as $\alpha(\textbf{x}, t) = \tilde{\alpha} (\textbf{x}) \exp[-i\omega t]$. Unlike a composite with a single matrix, the dissimilar thermal conductivity of a bimaterial leads to more rapid variation of the thermal field in the neighborhood of the bimaterial interface. Hence, it is necessary to employ interior interpolation points to improve the approximation accuracy for temperature distribution. Specifically, the numerical verifications have considered $128$ evenly distributed interior interpolation points in this paper. 

Using the reciprocal properties of bimaterial Green's function, two domain integrals in Eq. (\ref{eq:BIE_con}) can be converted into two surface integrals. Without the loss of any generality, this subsection only provides the derivation of domain integrals in $\mathcal{D}^+$, which can be extended to domain integrals in $\mathcal{D}^-$. Substituting Eq. (\ref{eq:backward_temp}) for time partial derivative approximation and Eq. (\ref{eq:interpolation_temp}) for temperature distribution approximation into the last two terms in Eq. (\ref{eq:BIE_con}), the domain integrals can be written as: 
\begin{small}
\begin{equation}
\begin{split}
    &\int_{\mathcal{D}^+} C_p' G(\textbf{x}, \textbf{x}') \frac{\partial u(\textbf{x}', t)}{\partial t} \thinspace d\textbf{x}'+ \int_{\mathcal{D}^-} C_p'' G(\textbf{x}, \textbf{x}') \frac{\partial u(\textbf{x}', t)}{\partial t} \thinspace d\textbf{x}'\\
    & = \frac{\sum_{m = 1}^{NN + NS} \left[ 3 \alpha(\textbf{x}^m, t_n) - 4 \alpha(\textbf{x}^m, t_{n-1}) + \alpha(\textbf{x}^m, t_{n-2}) \right] \Big\{ \int_{\mathcal{D}^+} C_p' \left[ G(\textbf{x}', \textbf{x}') \nabla^2 \Gamma(\textbf{x}, \textbf{x}^m) \right] \thinspace d\textbf{x}' + \int_{\mathcal{D}^-} C_p'' \left[ G(\textbf{x}', \textbf{x}') \nabla^2 \Gamma(\textbf{x}, \textbf{x}^m) \right] \thinspace d\textbf{x}' \Big\} }{2 \Delta t} 
\end{split}
    \label{eq:conv_dom}
\end{equation}
\end{small}
where the entire domain integrals have been written in terms of $NN + NS$ elementary integrals multiplied by the virtual source density at the corresponding interpolation point, and it only requires the storage of the virtual source density of two preceding time stations. Each elementary integral in Eq. (\ref{eq:conv_dom}) can be simplified through the reciprocal property of Green's function. For example, the domain integral in the upper material domain $\mathcal{D}^+$ can be written as:
\begin{equation}
\begin{split}
    & \int_{\mathcal{D}^+} C_p' G(\textbf{x}, \textbf{x}') \Gamma_{,i'i'}(\textbf{x}', \textbf{x}^m) \thinspace d\textbf{x}'  = \int_{\mathcal{D}} C_p' \frac{\partial}{\partial x_i'} \left[ G(\textbf{x}, \textbf{x}') \Gamma_{,i'}(\textbf{x}', \textbf{x}^m) \right] \thinspace d\textbf{x}' - \int_{\mathcal{D}} C_p' G_{,i'}(\textbf{x}, \textbf{x}') \Gamma_{,i'}(\textbf{x}', \textbf{x}^m) \thinspace d\textbf{x}' \\ 
    & = \int_{\partial \mathcal{D}^+} n_i(\textbf{x}') C_p' G(\textbf{x}, \textbf{x}') \Gamma_{,i'}(\textbf{x}', \textbf{x}^m) \thinspace d\textbf{x}' + \int_{\mathcal{D}} \frac{\partial}{\partial x_i'} C_p' \left\{  G_{,i'}(\textbf{x}, \textbf{x}') \Gamma(\textbf{x}', \textbf{x}^m) \right\} \thinspace d\textbf{x}' - \int_{\mathcal{D}^+} C_p' G_{,i'i'}(\textbf{x}, \textbf{x}') \Gamma(\textbf{x}', \textbf{x}^m) \thinspace d\textbf{x}' \\ 
    & = \int_{\partial \mathcal{D}^+} C_p' \Gamma_{,i'}(\textbf{x}', \textbf{x}^m) n_i(\textbf{x}') G(\textbf{x}, \textbf{x}') \thinspace d\textbf{x}' - \int_{\partial \mathcal{D}^+} C_p' G_{,i'}(\textbf{x}, \textbf{x}') n_i(\textbf{x}') \Gamma(\textbf{x}', \textbf{x}^m) \thinspace d\textbf{x}' -  c(\textbf{x}) \frac{C_p'}{K'} \Gamma(\textbf{x}, \textbf{x}^m)
\end{split}
    \label{eq:reci_again_up}
\end{equation}
Similarly, the lower material domain integral in $\mathcal{D}^-$ can be obtained. Therefore, the two domain integrals yield two surface integrals and two free terms:  
\begin{equation}
\begin{split}
     & \int_{\mathcal{D}^+} C_p' G(\textbf{x}, \textbf{x}') \Gamma_{,i'i'}(\textbf{x}', \textbf{x}^m) \thinspace d\textbf{x}' + \int_{\mathcal{D}^-} C_p'' G(\textbf{x}, \textbf{x}') \Gamma_{,i'i'}(\textbf{x}', \textbf{x}^m) \thinspace d\textbf{x}'  \\ & = \int_{\partial \mathcal{D}^+} C_p' \Gamma_{,i'}(\textbf{x}', \textbf{x}^m) n_i(\textbf{x}') G(\textbf{x}, \textbf{x}') \thinspace d\textbf{x}' - \int_{\partial \mathcal{D}^+} C_p' G_{,i'}(\textbf{x}, \textbf{x}') n_i(\textbf{x}') \Gamma(\textbf{x}', \textbf{x}^m) \thinspace d\textbf{x}' \\ & + \int_{\partial \mathcal{D}^-} C_p'' \Gamma_{,i'}(\textbf{x}', \textbf{x}^m) n_i(\textbf{x}') G(\textbf{x}, \textbf{x}') \thinspace d\textbf{x}' - \int_{\partial \mathcal{D}^-} C_p'' G_{,i'}(\textbf{x}, \textbf{x}') n_i(\textbf{x}') \Gamma(\textbf{x}', \textbf{x}^m) \thinspace d\textbf{x}' -  c(\textbf{x}) \frac{C_p(\textbf{x})}{K(\textbf{x})} \Gamma(\textbf{x}, \textbf{x}^m)
\end{split}
    \label{eq:reci_again}
\end{equation}
In Eq. (\ref{eq:reci_again}), the free term $c(\textbf{x})$ is combined together because the interior point can only exist in one of the bi-layered matrices. Specifically, when the interior field point is located in the upper phase, $c(\textbf{x}) \frac{C_p(\textbf{x})}{K(\textbf{x})} = C_p' / K'$, and vice versa. Substituting Eqs. (\ref{eq:reci_again}) and (\ref{eq:conv_dom}) into (\ref{eq:BIE_con}), the bimaterial dual-reciprocity boundary integral equations (DR-BIEs) for transient heat transfer can be formulated:
\begin{equation}
\begin{split}
    &c(\textbf{x}) u(\textbf{x}, t_n) + \int_{\partial \mathcal{D}} K(\textbf{x}') G_{,i'}(\textbf{x}, \textbf{x}') n_i(\textbf{x}') u^{BC}(\textbf{x}', t_n) \thinspace d\textbf{x}'  = \int_{\partial \mathcal{D}} G(\textbf{x}, \textbf{x}') q^{BC}(\textbf{x}', t_n) \thinspace d\textbf{x}' \\ & + \sum_{m = 1}^{NN + NS} \frac{1}{2 \Delta t} \left[ 3 \alpha(\textbf{x}^m, t_n) - 4 \alpha(\textbf{x}^m, t_{n-1}) + \alpha(\textbf{x}^m, t_{n-2}) \right] \Big\{ -\int_{\partial \mathcal{D}^+} C_p' \Gamma_{,i'}(\textbf{x}', \textbf{x}^m) n_i(\textbf{x}') G(\textbf{x}, \textbf{x}') \thinspace d\textbf{x}' \\ & + \int_{\partial \mathcal{D}^+} C_p' G_{,i'}(\textbf{x}, \textbf{x}') n_i(\textbf{x}') \Gamma(\textbf{x}', \textbf{x}^m) \thinspace d\textbf{x}' -\int_{\partial \mathcal{D}^-} C_p'' \Gamma_{,i'}(\textbf{x}', \textbf{x}^m) n_i(\textbf{x}') G(\textbf{x}, \textbf{x}') \thinspace d\textbf{x}' \\ &+ \int_{\partial \mathcal{D}^-} C_p'' G_{,i'}(\textbf{x}, \textbf{x}') n_i(\textbf{x}') \Gamma(\textbf{x}', \textbf{x}^m) \thinspace d\textbf{x}' + c(\textbf{x}) \frac{C_p(\textbf{x})}{K(\textbf{x})} \Gamma(\textbf{x}, \textbf{x}^m) \Big\}
\end{split}
    \label{eq:dual_BIE_transient}
\end{equation}

Following the same fashion, the substitution of Eqs. (\ref{eq:dif_harmonic}) and (\ref{eq:interpolation_temp}) into (\ref{eq:BIE_con}) yields the bimaterial DR-BIEs for the time-harmonic heat transfer as follows:
\begin{equation}
    \begin{split}
    &c(\textbf{x}) \tilde{u}(\textbf{x}) + \int_{\partial \mathcal{D}} K(\textbf{x}') G_{,i'}(\textbf{x}, \textbf{x}') n_i(\textbf{x}') \tilde{u}^{BC}(\textbf{x}') \thinspace d\textbf{x}'  = \int_{\partial \mathcal{D}} G(\textbf{x}, \textbf{x}') \tilde{q}^{BC}(\textbf{x}') \thinspace d\textbf{x}' \\ & -i \omega \sum_{m = 1}^{NN + NS} \frac{1}{2 \Delta t} \tilde{\alpha}(\textbf{x}^m) \Big\{ \int_{\partial \mathcal{D}} C_p(\textbf{x}') G_{,i'}(\textbf{x}, \textbf{x}') n_i(\textbf{x}') \Gamma(\textbf{x}', \textbf{x}^m) \thinspace d\textbf{x}' \\ & + \int_{\partial \mathcal{D}} C_p(\textbf{x}') G(\textbf{x}, \textbf{x}') \Gamma_{,i'}(\textbf{x}', \textbf{x}^m) n_i(\textbf{x}') \thinspace d\textbf{x}' + c(\textbf{x}) \frac{C_p(\textbf{x})}{K(\textbf{x}) } \Gamma(\textbf{x}, \textbf{x}^m) \Big\}
\end{split}
    \label{eq:dual_BIE_harmonic}
\end{equation}
where the transient and harmonic heat transfer problems employ the same boundary integrals, and the only difference lies in the time partial derivative of temperature. Note that using bimaterial Green's function in the boundary element method is not new in steady-state problems. For example, bi-layered elastic and thermoelastic analyses have been conducted in \cite{Wu2023_EABE} and \cite{Wu2024_IJES}, respectively. However, for transient or time-harmonic heat transfer, although the bimaterial DR-BIEs share two same surface integrals as previous publications, Eq. (\ref{eq:dual_BIE_transient}) or (\ref{eq:dual_BIE_harmonic}) involves time-/frequency-dependent source terms, which is new and requires a few interior nodes along the interface for improved accuracy considering the time-dependent discontinued temperature gradient. 

Comparing the bimaterial DR-BIEs with its homogeneous-material counterpart, the relative positions of source point $\textbf{x}'$ and field point \textbf{x} should be carefully taken into account for the surface integrals associated with the Green's function. Note that the bimaterial DR-BIEs are fundamentally different from the multi-region BIEs using the Green's function for an infinite domain as follows: (i) The multi-region scheme is built on continuity equations, which require the discretization of the interface with high resolution for accuracy. For instance, the rank of the global coefficient matrix for the multi-region scheme equals the number of nodes on the exterior surface, bimaterial interface, and interpolation points. (ii) The bimaterial DR-BIE scheme utilizes properties of Green's function, which sets virtual nodes/elements on the bimaterial interface to evaluate surface integrals. However, the nodal values can be approximated by the virtual source density and radial basis functions. In such a case, the bimaterial DR-BIE does not set any unknowns on the bimaterial interface. As a result, the rank of the global coefficient matrix for the bimaterial DR-BIE is determined by the nodes of the exterior surface and a few interior points, which is much less than that of the multi-region scheme at the equivalent accuracy, because the bimaterial Green's function addresses the interfacial continuity. 

\subsection{Bimaterial dual-reciprocity inclusion-based boundary element method}
This paper proposes to extend Eshelby's equivalent inclusion method (EIM) to the inhomogeneity problem of the bimaterials, in which each inhomogeneity is treated as an equivalent inclusion with the same material properties as the surrounding matrix, and the material mismatch is simulated by continuous eigen-fields. Specifically, the eigen-temperature-gradient (ETG) \cite{hatta1986} and the eigen-heat-source (EHS) \cite{Wu2024-prsa} are employed to simulate the mismatch of thermal conductivity and heat capacity, respectively. When the time interval $\Delta t$ is not large, it is rational to assume eigen-fields are constant between two adjacent time stations. Based on previous case studies on steady-state heat conduction \cite{Wang_2022_JAM}, when inhomogeneities are close to each other, or in the neighborhood of the bimaterial interface, ETGs have exhibited significant spatial variation due to intensive interactions. Moreover, our recent work \cite{Wu2024-prsa} has proved that for transient or time-harmonic heat transfer, eigen-fields cannot be uniform. Therefore, eigen-fields are defined as piece-wise polynomial functions within each subdomain: 
\begin{small}
\begin{equation}
\begin{aligned}
    u_i^*(\textbf{x}, t)  = &\sum_{n=1}^{N^t} \sum_{I=1}^{N^I} \left[ u_{i}^{I0*}(t_n) + u^{I1*}_{ip}(t_n) (x_p - x_p^{IC}) + u_{ipq}^{I2*}(t_n) (x_p - x_p^{IC}) (x_q - x_q^{IC}) \right]  \Theta(1 - \xi^I) H(t - t_{n-1}) H(t_n - t)  \\
    Q^*(\textbf{x}, t)  =& \sum_{n=1}^{N^t} \sum_{I=1}^{N^I} \left[ Q^{I0*}(t_n) + Q^{I1*}_{p}(t_n) (x_p - x_p^{IC}) + Q^{I2*}_{pq}(t_n) (x_p - x_p^{IC}) (x_q - x_q^{IC}) \right] \Theta(1 - \xi^I) H(t - t_{n-1}) H(t_n - t) 
\end{aligned}
    \label{eq:eigen-fields-poly}
\end{equation}
\end{small}
where $N^t$ is the number of time stations, and $H(t - t')$ is the Heaviside function; $u^{I0*}_{i}(t_n), u^{I1*}_{ip}(t_n), u^{I2*}_{ipq}(t_n)$ are uniform, linear, quadratic terms of the polynomial ETG expanded at the center of the $I^{th}$ subdomain at time station $t_n$; similarly, $Q^{I0*}(t_n), Q^{I1*}_{p}(t_n), Q^{I2*}_{pq}(t_n)$ are uniform, linear, quadratic terms of the polynomial EHS expanded at the center of the $I^{th}$ subdomain at $t_n$, and; $\Theta( 1- \xi)$ is the characteristic function of the $I^{th}$ ellipsoidal subdomain centered at the point $\textbf{x}^{IC}$: 
\begin{equation}
    \Theta(1 - \xi^I) = \begin{cases} 1 & \xi^I \in [0, 1) \\ 0 & \xi^I \in [1, +\infty)\end{cases} \quad; \quad \xi^I = \frac{(x_1 - x_1^{IC})^2}{(a_1^I)^2} + \frac{(x_2 - x_2^{IC})^2}{(a_2^I)^2} + \frac{(x_3 - x_3^{IC})^2}{(a_3^I)^2}
    \label{eq:Theta}
\end{equation}
For time-harmonic heat transfer, the spatial part of eigen-fields are written in a polynomial form as follows: 
\begin{equation}
\begin{aligned}
    \tilde{u}_i^*(\textbf{x}, t) = \sum_{I=1}^{N^I} \left[ \tilde{u}_{i}^{I0*} + \tilde{u}^{I1*}_{ip} (x_p - x_p^{IC}) + \tilde{u}_{ipq}^{I2*} (x_p - x_p^{IC}) (x_q - x_q^{IC}) \right] \Theta(1 - \xi^I)\\ 
    \tilde{Q}^*(\textbf{x}, t) = \sum_{I=1}^{N^I} \left[ \tilde{Q}^{I0*} + \tilde{Q}^{I1*}_{p}(x_p - x_p^{IC}) + \tilde{Q}^{I2*}_{pq} (x_p - x_p^{IC}) (x_q - x_q^{IC}) \right] \Theta(1 - \xi^I)\\ 
\end{aligned}
    \label{eq:eigen-harmonic}
\end{equation}
where $\tilde{u}^*$ and $\tilde{Q}^*$, and their polynomial coefficients are complex numbers to show the heat flux wave.

Based on the bimaterial Green's function, disturbances caused by the eigen-fields of ETG and EHS can be written in terms of domain integrals of Green's function multiplied by the source fields. At time station $t_n$, the disturbed temperature at the field point $\textbf{x}$ is:
\begin{equation}
    u'(\textbf{x}, t_n) = \sum_{I=1}^{N^I} \left[ \int_{\Omega^I} G_{,i'}(\textbf{x}, \textbf{x}') K(\textbf{x}') u_{i}^{I*}(\textbf{x}', t_n) d \thinspace \textbf{x}' -  \int_{\Omega^I}G(\textbf{x}, \textbf{x}') \frac{\partial Q^{I*}(\textbf{x}', t_n)}{\partial t} d \thinspace \textbf{x}' \right]
    \label{eq:disturbed_temp}
\end{equation}

As the domain integrals of bimaterial Green's function and polynomial-form source fields can be written in terms of bimaterial Eshelby's tensors \cite{Wu2021_jam_polygonal,Wu2021_polyhedral}, for transient  and time-harmonic heat transfer, the disturbed temperature can be, respectively, written as:
\begin{equation}
\begin{split}
    u'(\textbf{x}, t_n) = &\sum_{I=1}^{N^I} \Big\{ \left[ D_i u_i^{I0*}(t_n) + D_{ip} u_{ip}^{I1*}(t_n) + D_{ipq} u_{ipq}^{I2*}(t_n) \right] + \frac{L^I}{2 \Delta t} \left[ 3 Q^{I0*}(t_n) - 4 Q^{I0*}(t_{n-1}) + Q^{I0*}(t_{n-2}) \right] \\ & + \frac{L_{p}^I}{2 \Delta t} \left[ 3 Q_{p}^{I1*}(t_n) - 4 Q_{p}^{I1*}(t_{n-1}) + Q_{p}^{I1*}(t_{n-2}) \right] + \frac{L_{pq}^I}{2 \Delta t} \left[ 3 Q_{pq}^{I2*}(t_n) - 4 Q_{pq}^{I2*}(t_{n-1}) + Q_{pq}^{I2*}(t_{n-2}) \right] \Big\} 
\end{split}
    \label{eq:disturbed_temp_transient}
\end{equation}
and 
\begin{equation}
    \tilde{u}'(\textbf{x}) = \sum_{I=1}^{N^I} \Big\{ \left( D_i \tilde{u}_i^{I0*} + D_{ip} \tilde{u}_{ip}^{I1*} + D_{ipq} \tilde{u}_{ipq}^{I2*} \right) + \left( L^I Q^{I0*} + L^I_{p} Q_{p}^{I1*} + L^{I} Q_{pq}^{I2*} \right) \Big\}
    \label{eq:disturbed_temp_harmonic}
\end{equation}
where Eshelby's tensors are defined as:
\begin{equation}
\begin{aligned}
    & D_{ipq...}^I = \int_{\Omega^I} K(\textbf{x}') G_{,i'}(\textbf{x}, \textbf{x}') (x_p' - x_p^{IC}) (x_q' - x_q^{IC})... \thinspace d\textbf{x}' = \begin{cases} \frac{-M_I}{4 \pi} \left( \Phi_{pq...,i} + M_P M_Q... \frac{K^s - \overline{K}^s}{K^s + \overline{K}^s} \overline{\Phi}_{pq...,i} \right) & x_3 x_3^{IC} \geq 0 \\ \frac{-K^s}{2\pi (K^s + \overline{K}^s)} \Phi_{pq...,i} & x_3 x_3^{IC} < 0 \end{cases} \\
    & L_{pq...}^I = \int_{\Omega^I} G(\textbf{x}') (x_p' - x_p^{IC}) (x_q' - x_q^{IC})... \thinspace d\textbf{x}' = \begin{cases} \frac{1}{4\pi K^s} \left( \Phi_{pq...} + M_P M_Q... \frac{K^s - \overline{K}^s}{K^s + \overline{K}^s} \overline{\Phi}_{pq...} \right) & x_3 x_3^{IC} \geq 0 \\ \frac{1}{2 (K^s + \overline{K}^s)} \Phi_{pq...} & x_3 x_3^{IC} < 0 \end{cases}
\end{aligned}
    \label{eq:Eshelby_tensors}
\end{equation}
in which $D^I_{ipq...}$ and $L^I_{pq...}$ are the polynomial-form bimaterial Eshelby's tensors for ETG and EHS in the $I^{th}$ subdomain, respectively; $\Phi_{pq...} = \int_{\Omega^I} \phi (x_p' - x_p^{IC}) (x_q' - x_q^{IC})... \thinspace d\textbf{x}'$ and $\overline{\Phi}_{pq...} = \int_{\Omega^I} \overline{\phi} (x_p' - x_p^{IC}) (x_q' - x_q^{IC})... \thinspace d\textbf{x}'$ are the domain integrals of $\phi$ and $\overline{\phi}$ multiplied with source densities \cite{Wang_2022_JAM}, and; the capital characters of indices, i.e., $M_P$ or $M_Q$, are not involved in the dummy index rule; definitions of superscript $s$ are consistent with Eq. (\ref{eq:Gfunc}). Note that Eq. (\ref{eq:disturbed_temp_transient}) shares the same formulae as the corresponding expressions of a single-material. However, each Eshelby's tensor should be replaced with the bimaterial case in Eq. (\ref{eq:Eshelby_tensors}). Unlike Eshelby's tensors evaluated with Green's function for a single material, the bimaterial Eshelby's tensors are functions of three semi-axes of the ellipsoidal subdomain and its relative position to the bimaterial interface. 

Superposing the influence from boundary responses and disturbances by eigen-fields, the temperature at time station $t_n$ and field point $\textbf{x}$ can be written as the combination of Eqs. (\ref{eq:dual_BIE_transient}) and  (\ref{eq:disturbed_temp_transient}) for transient heat transfer:
\begin{equation}
\begin{split}
    &c(\textbf{x}) u(\textbf{x}, t_n) + \int_{\partial \mathcal{D}} K(\textbf{x}') G_{,i'}(\textbf{x}, \textbf{x}') n_i(\textbf{x}') u^{BC}(\textbf{x}', t_n) \thinspace d\textbf{x}'  = \int_{\partial \mathcal{D}} G(\textbf{x}, \textbf{x}') q^{BC}(\textbf{x}', t_n) \thinspace d\textbf{x}' \\ & + \sum_{m = 1}^{NN + NS} \frac{1}{2 \Delta t} \left[ 3 \alpha(\textbf{x}^m, t_n) - 4 \alpha(\textbf{x}^m, t_{n-1}) + \alpha(\textbf{x}^m, t_{n-2}) \right] \Big\{ -\int_{\partial \mathcal{D}^+} C_p' \Gamma_{,i'}(\textbf{x}', \textbf{x}^m) n_i(\textbf{x}') G(\textbf{x}, \textbf{x}') \thinspace d\textbf{x}' \\ & + \int_{\partial \mathcal{D}^+} C_p' G_{,i'}(\textbf{x}, \textbf{x}') n_i(\textbf{x}') \Gamma(\textbf{x}', \textbf{x}^m) \thinspace d\textbf{x}' -\int_{\partial \mathcal{D}^-} C_p'' \Gamma_{,i'}(\textbf{x}', \textbf{x}^m) n_i(\textbf{x}') G(\textbf{x}, \textbf{x}') \thinspace d\textbf{x}' \\ &+ \int_{\partial \mathcal{D}^-} C_p'' G_{,i'}(\textbf{x}, \textbf{x}') n_i(\textbf{x}') \Gamma(\textbf{x}', \textbf{x}^m) \thinspace d\textbf{x}' + c(\textbf{x}) \frac{C_p(\textbf{x})}{K(\textbf{x})} \Gamma(\textbf{x}, \textbf{x}^m) \Big\} \\ & + \sum_{I=1}^{N^I} \Big\{ \left[ D_i u_i^{I0*}(t_n) + D_{ip} u_{ip}^{I1*}(t_n) + D_{ipq} u_{ipq}^{I2*}(t_n) \right] + \frac{L^I}{2 \Delta t} \left[ 3 Q^{I0*}(t_n) - 4 Q^{I0*}(t_{n-1}) + Q^{I0*}(t_{n-2}) \right] \\ & + \frac{L_{p}^I}{2 \Delta t} \left[ 3 Q_{p}^{I1*}(t_n) - 4 Q_{p}^{I1*}(t_{n-1}) + Q_{p}^{I1*}(t_{n-2}) \right] + \frac{L_{pq}^I}{2 \Delta t} \left[ 3 Q_{pq}^{I2*}(t_n) - 4 Q_{pq}^{I2*}(t_{n-1}) + Q_{pq}^{I2*}(t_{n-2}) \right] \Big\}
\end{split}
    \label{eq:overll_transient}
\end{equation}

Following the same fashion, the combination of Eqs. (\ref{eq:dual_BIE_harmonic}) and (\ref{eq:disturbed_temp_transient}) yield the spatial variation of temperature at field point $\textbf{x}$ for time-harmonic heat transfer:

\begin{equation}
\begin{split}
    &c(\textbf{x}) \tilde{u}(\textbf{x}) + \int_{\partial \mathcal{D}} K(\textbf{x}') G_{,i'}(\textbf{x}, \textbf{x}') n_i(\textbf{x}') \tilde{u}^{BC}(\textbf{x}') \thinspace d\textbf{x}'  = \int_{\partial \mathcal{D}} G(\textbf{x}, \textbf{x}') \tilde{q}^{BC}(\textbf{x}') \thinspace d\textbf{x}' \\ & -i \omega \sum_{m = 1}^{NN + NS} \frac{1}{2 \Delta t} \tilde{\alpha}(\textbf{x}^m) \Big\{ \int_{\partial \mathcal{D}} C_p(\textbf{x}') G_{,i'}(\textbf{x}, \textbf{x}') n_i(\textbf{x}') \Gamma(\textbf{x}', \textbf{x}^m) \thinspace d\textbf{x}' \\ & + \int_{\partial \mathcal{D}} C_p(\textbf{x}') G(\textbf{x}, \textbf{x}') \Gamma_{,i'}(\textbf{x}', \textbf{x}^m) n_i(\textbf{x}') \thinspace d\textbf{x}' + c(\textbf{x}) \frac{C_p(\textbf{x})}{K(\textbf{x}) } \Gamma(\textbf{x}, \textbf{x}^m) \Big\}  \\ & + \sum_{I=1}^{N^I} \Big\{ \left( D_i \tilde{u}_i^{I0*} + D_{ip} \tilde{u}_{ip}^{I1*} + D_{ipq} \tilde{u}_{ipq}^{I2*} \right) + \left( L^I Q^{I0*} + L^I_{p} Q_{p}^{I1*} + L^{I}_{pq} Q_{pq}^{I2*} \right) \Big\}
\end{split}
    \label{eq:overll_harmonic}
\end{equation}

As thermal fields can be expressed by the bimaterial DR-BIEs and disturbances by eigen-fields within inclusions, the ETG is determined by equivalent flux conditions \cite{hatta1986}, which are constructed at the center of each subdomain. 
\begin{equation}
\begin{aligned}
    -K^s \left[ u_{,i}(\textbf{x}^{IC}, t_n) - u^{I0*}_{i}(t_n) \right] & = -K^I u_{,i}(\textbf{x}^{IC}, t_n) \\
    -K^s \left[ u_{,ip}(\textbf{x}^{IC}, t_n) - u^{I1*}_{ip}(t_n) \right] & = -K^I u_{,ip}(\textbf{x}^{IC}, t_n)\\ 
    -K^s \left[ u_{,ipq}(\textbf{x}^{IC}, t_n) - 2 u^{I2*}_{ipq}(t_n) \right] & = -K^I u_{,ipq}(\textbf{x}^{IC}, t_n)
\end{aligned}
\label{eq:equiv_flux_transient}
\end{equation}

Similarly, the EHS can be determined by equivalent heat generation/storage conditions \cite{Wu2024-prsa} as follows:
\begin{equation}
\begin{aligned}
    C_p^s \left[ u(\textbf{x}^{IC}, t_n) - Q^{I0*}(t_n) \right] & = C_p^I u_{}(\textbf{x}^{IC}, t_n) \\
    C_p^s \left[ u_{,p}(\textbf{x}^{IC}, t_n) - Q^{I1*}_{ip}(t_n) \right] & = C_p^I u_{,p}(\textbf{x}^{IC}, t_n) \\ 
    C_p^s \left[ u_{,pq}(\textbf{x}^{IC}, t_n) - 2 Q^{I2*}_{ipq}(t_n) \right] & = C_p^I u_{,pq}(\textbf{x}^{IC}, t_n) 
\end{aligned}
    \label{eq:equiv_heat_transient}
\end{equation}
where the superscript $s = '$ for $x_3^{IC} \geq 0$ and $s = ''$ for $x_3^{IC} < 0$.  Therefore, the original heterogeneous bi-layered composite system has been handled by $NN$ boundary nodal unknowns, $NS$ interior interpolation points, and $N^I \times num$ unknown coefficients of eigen-fields, where $num = 4, 16, 40$ for uniform, linear, and quadratic assumption of eigen-fields, respectively. 

Similarly, for time harmonic heat transfer, the equivalent inclusion conditions can be written as
\begin{equation}
\begin{aligned}
     \quad -K^s \left[ \tilde{u}_{,i}(\textbf{x}^{IC}) - \tilde{u}^{I0*}_{i} \right]  = -K^I \tilde{u}_{,i}(\textbf{x}^{IC}); \quad & C_p^s \left[ \tilde{u}(\textbf{x}^{IC}) - \tilde{Q}^{I0*} \right] = C_p^I \tilde{u}(\textbf{x}^{IC}) \\
    -K^s \left[ \tilde{u}_{,ip}(\textbf{x}^{IC}) - \tilde{u}^{I1*}_{ip} \right] = -K^I \tilde{u}_{,ip}(\textbf{x}^{IC}); \quad & C_p^s \left[ \tilde{u}_{,p}(\textbf{x}^{IC}) - \tilde{Q}^{I1*}_{ip} \right] = C_p^I \tilde{u}_{,p}(\textbf{x}^{IC})\\ 
    -K^s \left[ \tilde{u}_{,ipq}(\textbf{x}^{IC}) - 2 \tilde{u}^{I2*}_{ipq} \right] = -K^I \tilde{u}_{,ipq}(\textbf{x}^{IC}); \quad & C_p^s \left[ \tilde{u}_{,pq}(\textbf{x}^{IC}) - 2 \tilde{Q}^{I2*}_{ipq} \right] = C_p^I \tilde{u}_{,pq}(\textbf{x}^{IC})
\end{aligned}
\label{eq:equiv_Inclusion_Harmonic}
\end{equation}

\section{Numerical verification}
This section aims to provide numerical verification on (i) bimaterial dual-reciprocity boundary integral equations (DR-BIEs) without inhomogeneities and (ii) the dual-reciprocity inclusion-based boundary element method (DR-iBEM) proposed in the previous section. Note that the transient heat transfer and time-harmonic heat transfer share the same boundary integral equations with the same boundary coefficient matrices, bimaterial Eshelby's tensors, and their only difference lies in the evaluation of time partial derivatives of temperature, which leads to slightly different arrangement of global matrices. Therefore, this section only shows verification and comparison of transient heat transfer, as the time-harmonic heat transfer is a special case of it. Without loss of any generality, thermal properties of the bimaterial and inhomogeneities are selected as: (i) the upper layer: $K' = 4 W / m \cdotp K$, $C_p' = 10 W / m^3 \cdotp K$; (ii) the lower layer: $K'' = 2 W / m \cdotp K$, $C_p'' = 3 W / m^3 \cdotp K$; and (iii) inhomogeneities: $K^I = 10 W / m \cdotp K$, $C_p^I = 1 W / m^3 \cdotp K$. As shown in Fig. \ref{fig:geometryc}, the dimensions $l_a = l_b = 1$ m, $h_1 = h_2 = 1$ m are assigned. The boundary conditions are set up as: (i) a sinusoidal temperature function $u^{BC}(t) = 10 \sin [\frac{\pi}{10} t]$ is uniformly applied on the top surface; (ii) zero temperature is prescribed on the bottom surface; and (iii) all other surfaces are adiabatic and initial temperature is zero. For two case studies, the bimaterial DR-BIEs / DR-iBEM utilized the same boundary mesh that $1,002$ boundary nodes, $1,000$ quadrilateral elements, and $128$ uniformly distributed interior interpolation points. 

\subsection{Verification of the DR-BIEs for transient heat transfer}
For the purpose of reproduction, the finite element method (FEM) adopted a uniform discretization of the bi-layered matrix, which used $332,045$ nodes and $78,608$ elements. Due to a large number of unknowns, the FEM consumes $6.85$ GB of RAM and takes $868$ seconds for simulation, while the DR-BIEs only occupy $0.74$ GB of RAM and take $73$ seconds for simulation. 

\begin{figure}
\begin{subfigure}{.5\textwidth}
\centering
\includegraphics[width = 1\linewidth,height = \textheight,keepaspectratio]{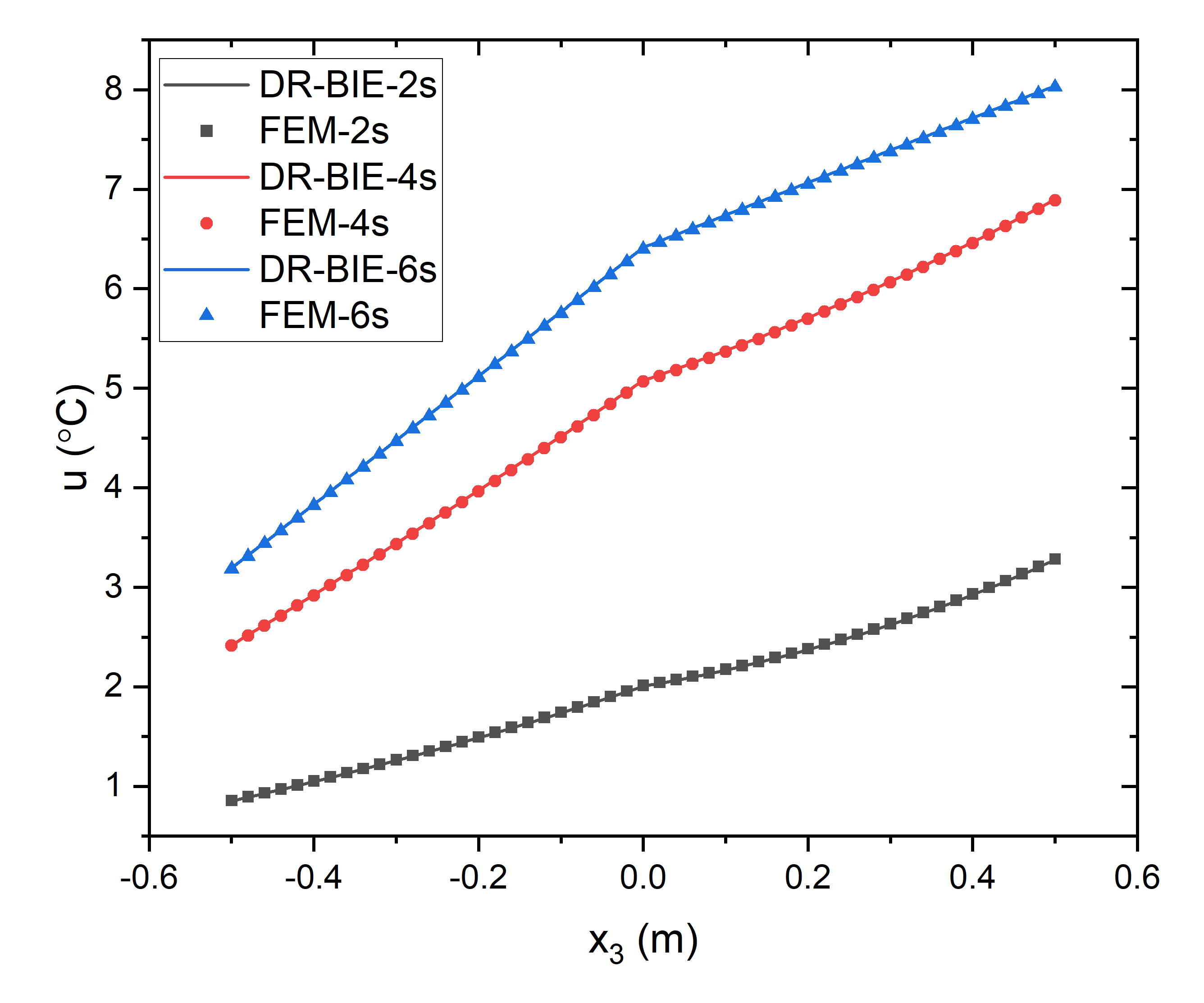}
\caption{}
\end{subfigure}
~
\begin{subfigure}{.5\textwidth}
\centering
\includegraphics[width = 1\linewidth,height = \textheight,keepaspectratio]{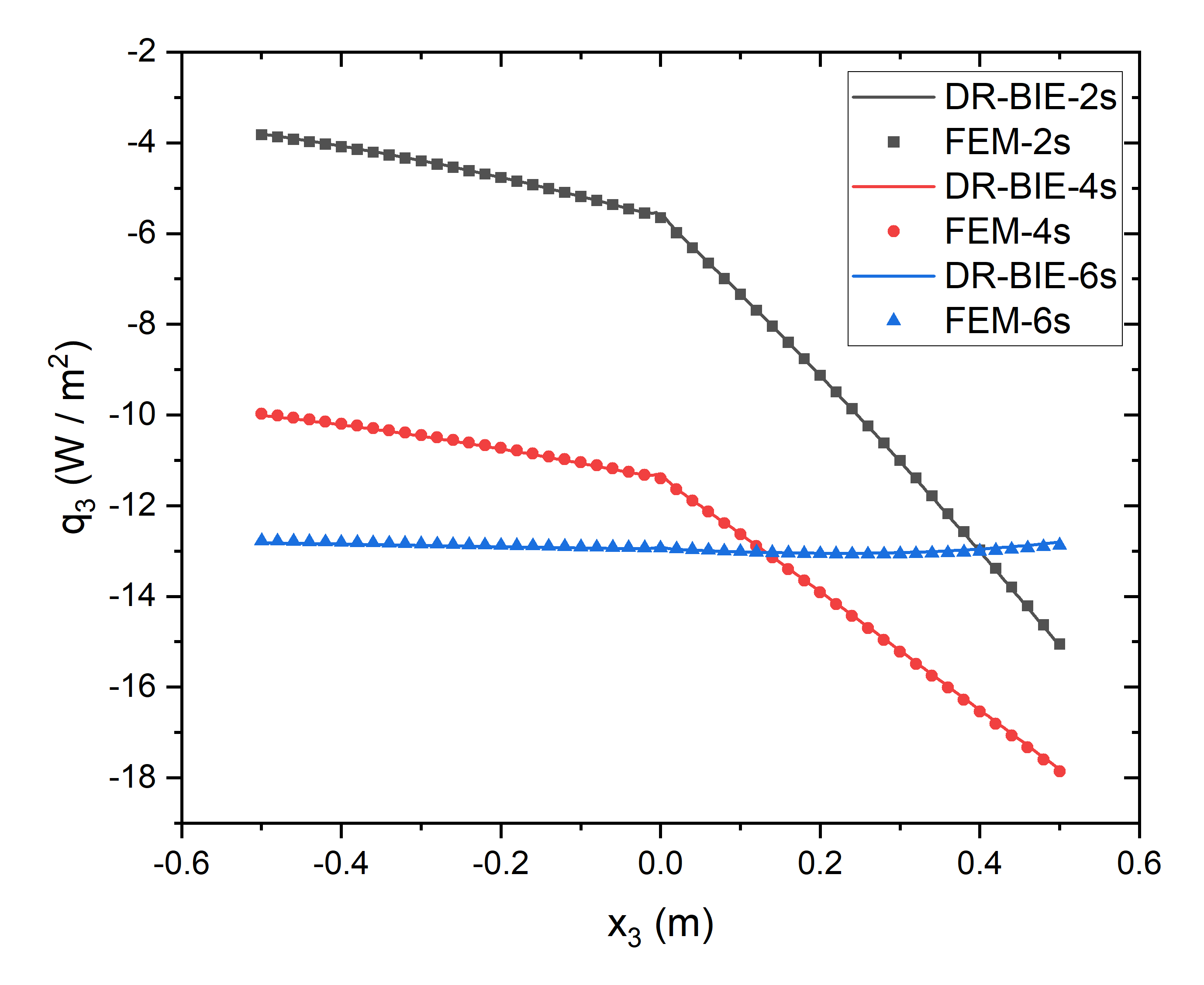}
\caption{}
\end{subfigure}
\caption{Variation and comparison of (a) temperature and (b) heat flux along the centerline $x_3 \in [-0.5, 0.5]$ m evaluated by the bimaterial dual-reciprocity boundary integral equations (DR-BIEs) and the finite element method (FEM) at time $t = 2, 4, 6$ s.}
\label{fig:thermal_bie}
\end{figure}

Figs. \ref{fig:thermal_bie} (a-b) show spatial variations of temperature and heat flux along the vertical centerline $x_3 \in [-0.5, 0.5]$ m at time $t = 2,4,6$ s, respectively. Fig. \ref{fig:thermal_bie} (a) shows that the temperature curves at time $t = 2, 4, 6$ s all exhibit continuous distribution across the bimaterial interface ($x_3 = 0$), which satisfy the continuity condition on the temperature in Section 2. Because the upper and lower layers exhibit different thermal conductivity, the temperature gradient in the third direction is not continuous at the bimaterial interface, which has been reported in steady-state heat conduction between layered materials. However, unlike steady-state heat conduction, temperature gradients change with time accordingly. For instance, when time $t = 2$ s, the temperature gradients in two phases are less than that of $4$ s and $6$ s. Note that the temperature on the top surface is a sine function, which increases from 0 to its maximum at 5 s (quarter period). This explains why the temperature gradients at time $2$ s are less than those at time $4$ s. However, the temperature gradients at 6 s are greater than those at 4 s because the heat transfer process depends on both thermal diffusivity and time duration, which leads to a temperature lag along the height direction. 

Fig. \ref{fig:thermal_bie} (b) shows that the third component of heat flux $q_3$ is continuous across the bimaterial interface, which satisfies the continuity condition on heat flux in Section 2. The heat flux continuity ensures that the temperature gradients on the upper and lower side of the bimaterial interface exhibit a ratio at the reverse of $K' / K''$, which is consistent with discontinuous temperature gradients in Fig. \ref{fig:thermal_bie} (a). 

\begin{figure}
\begin{subfigure}{.5\textwidth}
\centering
\includegraphics[width = 1\linewidth,height = \textheight,keepaspectratio]{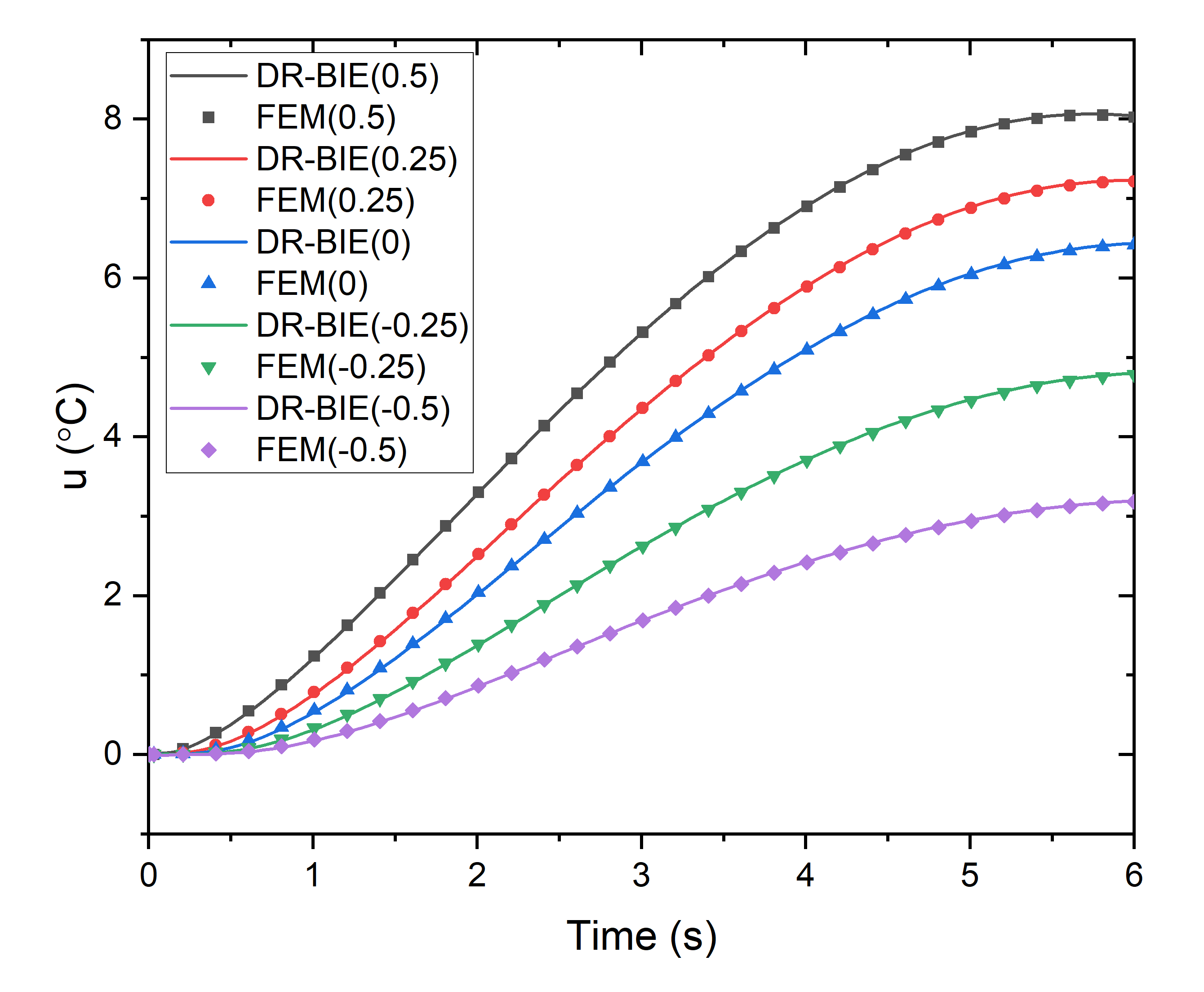}
\caption{}
\end{subfigure}
~
\begin{subfigure}{.5\textwidth}
\centering
\includegraphics[width = 1\linewidth,height = \textheight,keepaspectratio]{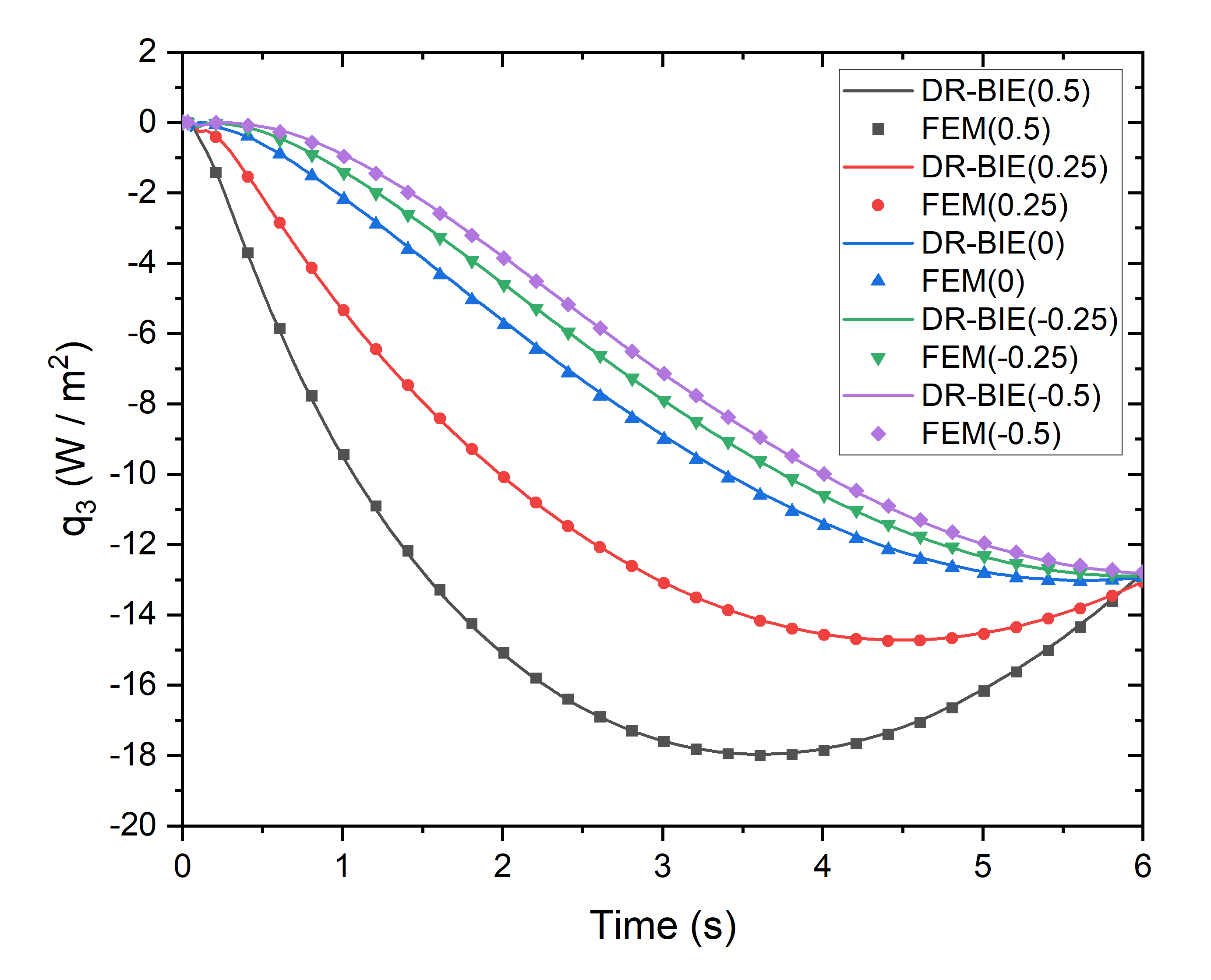}
\caption{}
\end{subfigure}
\caption{Variation and comparison of (a) temperature and (b) heat flux at field point (0.5, 0.5, $x_3$) ($x_3 = 0.5, 0.25, 0, -0.25, -0.5$ m) evaluated by the bimaterial dual-reciprocity boundary integral equations (DR-BIEs) and the finite element method (FEM), when time $t \in [0, 6]$ s.}
\label{fig:thermal_bie_time}
\end{figure}

Figs. \ref{fig:thermal_bie_time} (a-b) compare the temperature and heat flux variation with time at the height $x_3 = 0.5, 0.25, 0, -0.25, 0.5$ m when time $t \in [0, 6]$ s, respectively. Very minor discrepancies exist among results obtained by DR-BIEs and FEM, with the maximum difference occuring at $5$ s on the two curves ``DR-BIE (0.5)'' in Figs. \ref{fig:thermal_bie_time} (a-b), which is less than $0.1\%$. The temperature lag can be observed in Figs. \ref{fig:thermal_bie_time} (a-b) that temperature at greater heights exhibit larger magnitude and takes less time to reach its maxima. For instance, the curve ``DR-BIE(0.5)'' in Fig. \ref{fig:thermal_bie_time} (b) reaches the maxima approximately 1 s earlier than the curve ``DR-BIE(0.25)''. Considering the RAM usage and time consumption, the bimaterial DR-BIEs exhibit much better simulation accuracy and computational efficiency compared to the FEM in bi-layered transient heat transfer. It should be noted that the current numerical case study only considers a homogeneous bi-layered sample, and thus the FEM employs a uniform mesh strategy. However, the existence of inhomogeneities will cause issues with mesh transition/convergence, which inevitably leads to a larger number of degrees of freedom.

\subsection{Verification of the DR-iBEM for transient heat transfer}
This subsection compares the DR-iBEM and FEM on ellipsoidal inhomogeneity problems. As demonstrated in the previous subsection, dissimilar thermal properties play a crucial role in the heat transfer process across the bimaterial interface, resulting in sudden changes in temperature gradients. Therefore, to examine the accuracy of the DR-iBEM under intensive bimaterial interfacial effects and inhomogeneities' interactions, two cases are considered: (i) two equal-sized spherical inhomogeneities with radius $\textbf{a}^1 = \textbf{a}^2 = (0.1, 0.1, 0.1)$ m; (ii) two equal-sized spheroidal inhomogeneities with three semi-axes $\textbf{a}^1 = \textbf{a}^2 = (0.2, 0.2, 0.1)$ m are placed in the neighborhood of the bimaterial interface that: $\textbf{x}^{1C} = (0.5, 0.5, 0.125)$ and $\textbf{x}^{2C} = (0.5, 0.5, -0.125)$ m. 

The FEM simulation details are provided for two cases: (i) $980,759$ nodes, $715,454$ elements, and $21.74$ Gb RAM, which took $3,760$ s; (ii) $814,667$ nodes and $599,127$ elements, and $14.27$ Gb RAM, which took $2,895$ s. In contrast, the DR-iBEM with the quadratic order of eigen-fields used: (i) $0.87$ Gb for $132$ s; (ii) $0.87$ Gb for $137$ s,  including all post-process time duration and saved matrices to expedite the post-process . Note that the DR-iBEM adopts the same boundary discretization as the previous subsection, and no mesh is required for inhomogeneities, as they are simulated using eigen-fields. The solution time in case (ii) is slightly longer than that of the case (i) due to the spheroidal integrals, which involves the evaluation of the parameter $\lambda$ and its associated functions \cite{Moschovidis1975, Yang2025}.

\begin{figure}
\begin{subfigure}{.5\textwidth}
\centering
\includegraphics[width = 1\linewidth,height = \textheight,keepaspectratio]{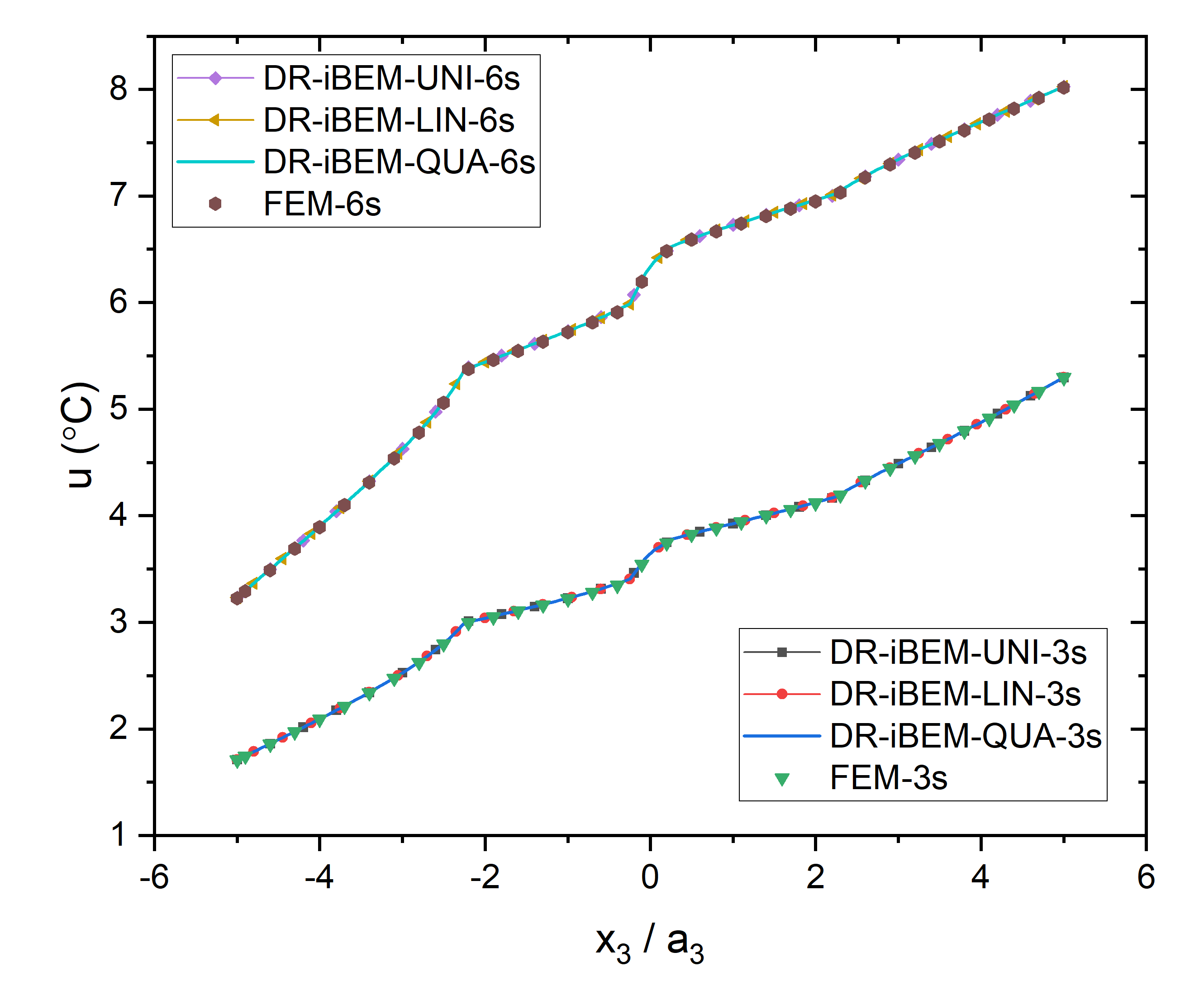}
\caption{}
\end{subfigure}
~
\begin{subfigure}{.5\textwidth}
\centering
\includegraphics[width = 1\linewidth,height = \textheight,keepaspectratio]{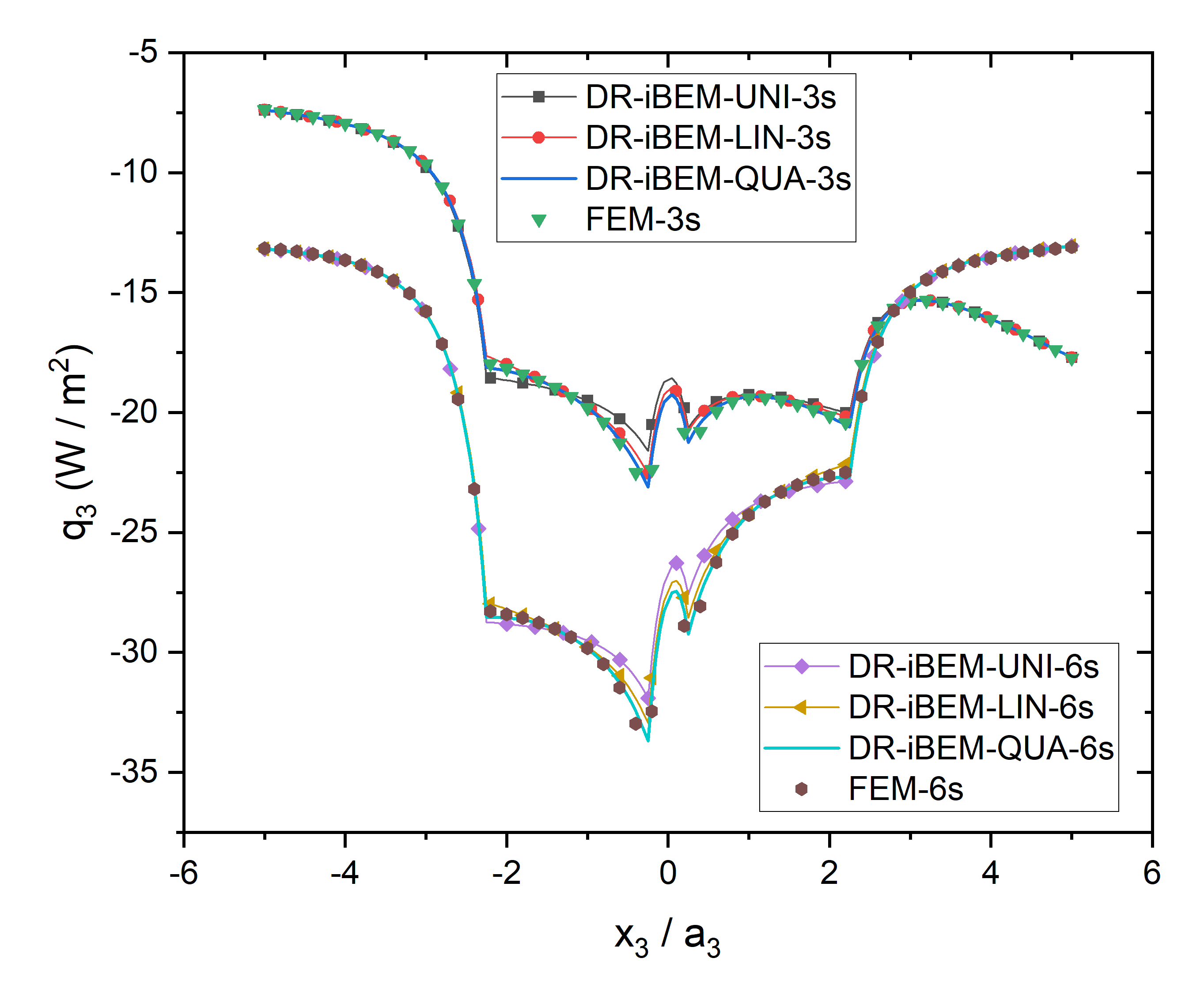}
\caption{}
\end{subfigure}
\caption{Variation and comparison of (a) temperature and (b) heat flux caused by two spherical inhomogeneities ($\textbf{a}^1 = \textbf{a}^2 = (0.1, 0.1, 0.1)$ m) along the centerline $x_3 \in [-0.5, 0.5]$ m evaluated by the dual-reciprocity inclusion-based boundary element method (DR-iBEM) with three order of eigen-fields and the finite element method (FEM) at time $t = 3, 6$ s.}
\label{fig:thermal_spheres}
\end{figure}

\begin{figure}
\begin{subfigure}{.5\textwidth}
\centering
\includegraphics[width = 1\linewidth,height = \textheight,keepaspectratio]{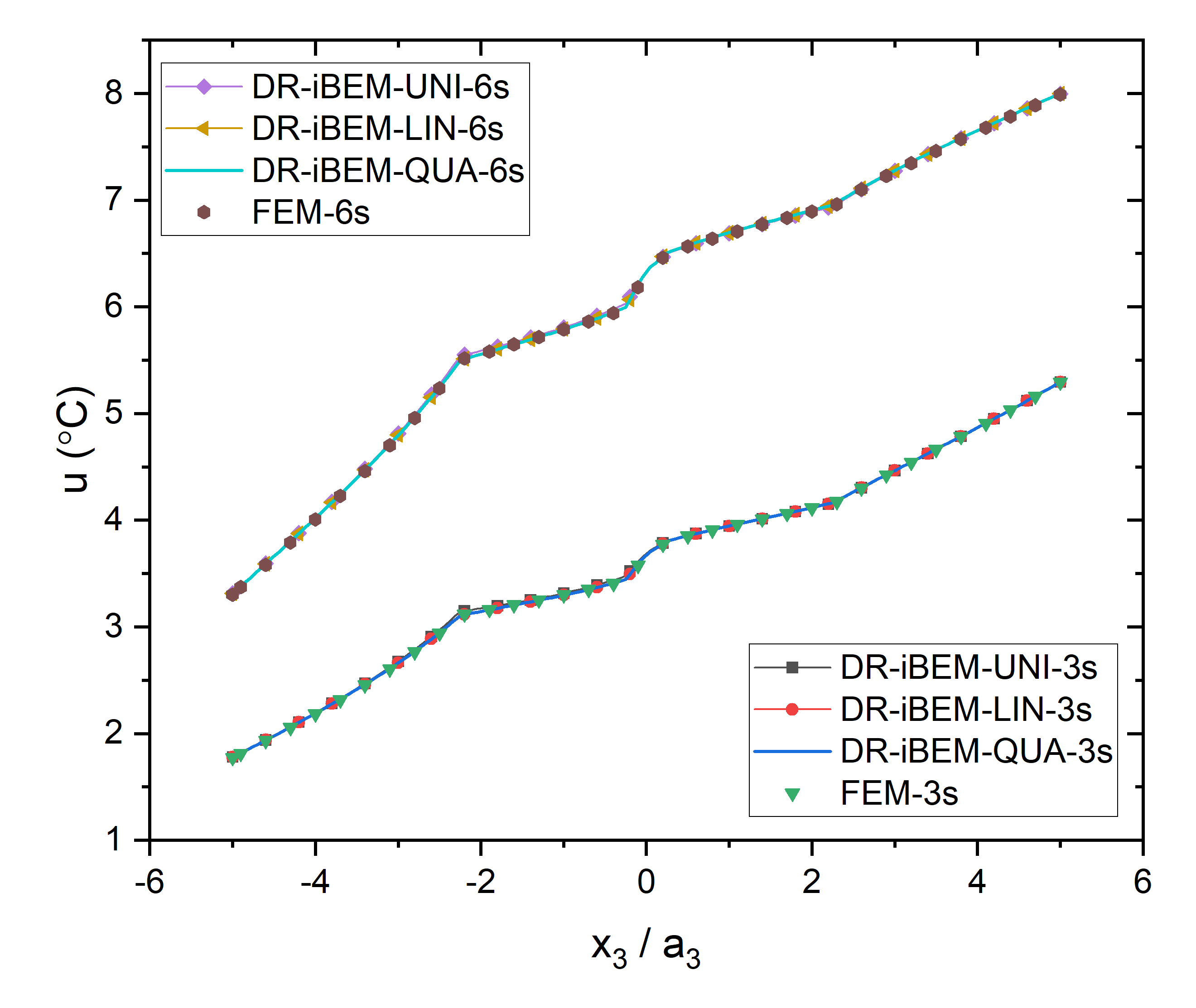}
\caption{}
\end{subfigure}
~
\begin{subfigure}{.5\textwidth}
\centering
\includegraphics[width = 1\linewidth,height = \textheight,keepaspectratio]{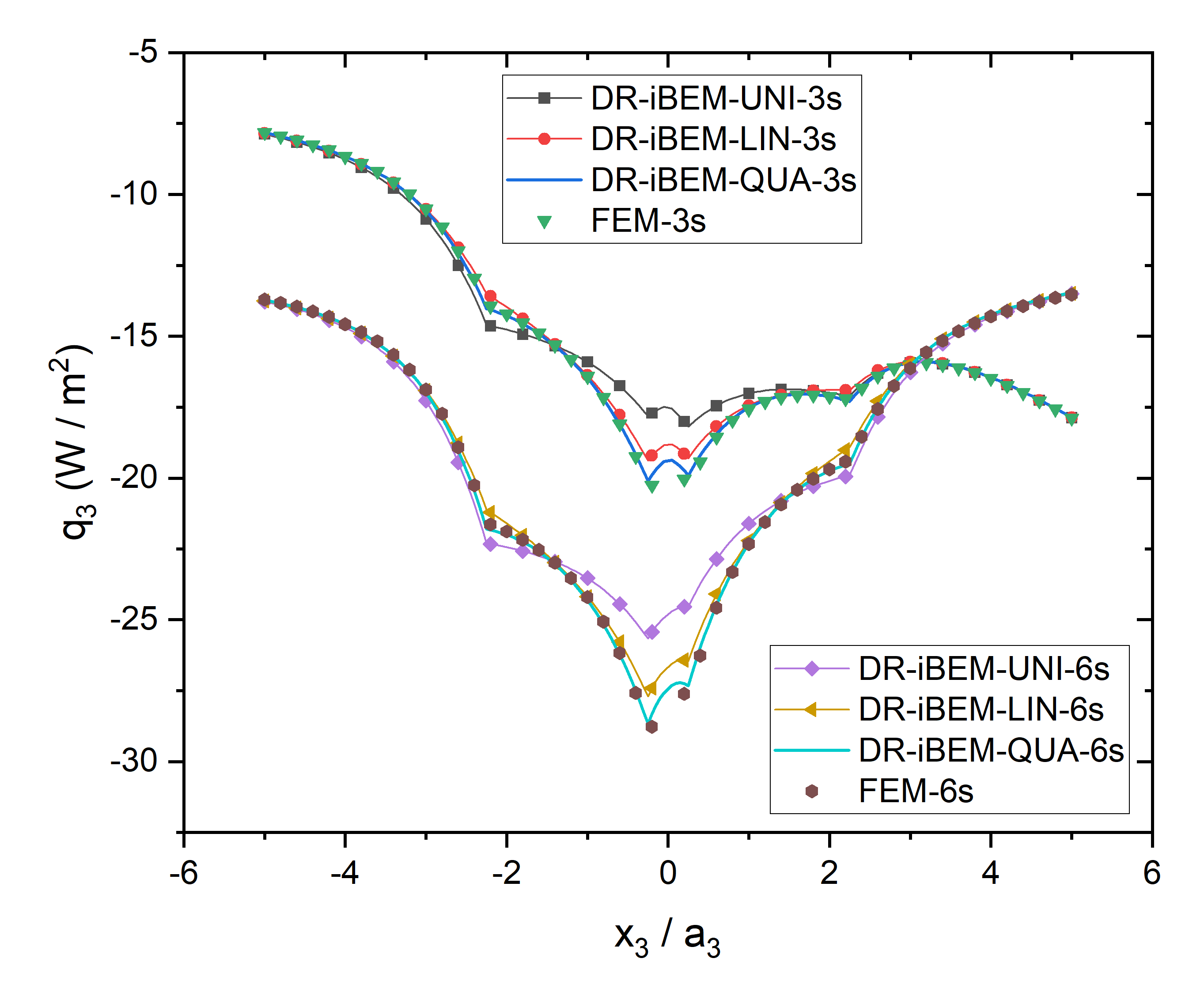}
\caption{}
\end{subfigure}
\caption{Variation and comparison of (a) temperature and (b) heat flux caused by two spheroidal inhomogeneities ($\textbf{a}^1 = \textbf{a}^2 = (0.2, 0.2, 0.1)$ m) along the centerline $x_3 \in [-0.5, 0.5]$ m evaluated by the dual-reciprocity inclusion-based boundary element method (DR-iBEM) with three order of eigen-fields and the finite element method (FEM) at time $t = 3, 6$ s.}
\label{fig:thermal_spheroids}
\end{figure}

Figs. \ref{fig:thermal_spheres} (a) and \ref{fig:thermal_spheroids} (a) compare the variation of temperature among DR-iBEM with three order eigen-fields and FEM at time $3, 6$ s, respectively. Generally, the differences between the DR-iBEM and FEM results are negligible, However, minor discrepancies around the entering region of two inhomogeneities can be found among two curves ``DR-iBEM-UNI-3s'', ``DR-iBEM-UNI-6s'', ``FEM-3s'', and ``FEM-6s'' with the errors less than $0.1\%$. Therefore, even under intensive bimaterial interfacial effects and inhomogeneities' interactions, the assumption of uniform eigen-fields can provide acceptable predictions on temperature fields. Since inhomogeneities exhibit greater thermal conductivity, the temperature gradients within inhomogeneities are less than those of the bi-layered matrix. Note that the temperature gradients on the bimaterial interface are discontinuous, which is caused by the different thermal conductivities in the upper and lower layers of materials. 

Figs. \ref{fig:thermal_spheres} (b) and \ref{fig:thermal_spheroids} (b) compare the variation of heat flux among DR-iBEM with three order eigen-fields and FEM at time $3, 6$ s, respectively. Although the uniform eigen-fields can provide acceptable predictions on temperature, obvious discrepancies between DR-iBEM-UNI and FEM illustrate the necessity to consider higher-order terms to describe the spatial variation of eigen-fields. In general, DR-iBEM with linear and quadratic eigen-fields can provide good predictions in most regions, except the entering points of the inhomogeneity, which is due to intensive interactions. Similar patterns have been reported in \cite{Wang_2022_JAM} for steady-state heat transfer. The discussion on the variation of eigen-fields began in \cite{Moschovidis1975}, in which the authors demonstrated that intensive inhomogeneities' interactions can significantly disturb the eigenstrain in elastostatics. Following their work, Wu et al. \cite{Wu2024_IJES} revealed that ETG distribution for a single inhomogeneity in a bimaterial is not uniform anymore due to the boundary/interface effects and particle-boundary interactions. By adding another inhomogeneity in the neighborhood, the spatial variation of eigen-fields is more significant, which is similar to the particle interactions in elastic wave propagation. Wu et al. \cite{Wu2024-prsa} proved that there does not exist any uniform far-field temperature gradients in unsteady heat transfer, and thus, a uniform ETG cannot satisfy the equivalent flux conditions at interior points of an ellipsoidal inhomogeneity. Furthermore, although the DR-iBEM utilizes closed-form Eshelby's tensors, the time convolution effects of eigen-fields are implicitly considered by the finite difference scheme, which treats two previous steps as pseudo-initial conditions.

\section{Application to a ceramic-metal FGM made of alumina and nickel}
As the DR-iBEM algorithm has been verified to address inhomogeneity problems with high accuracy and computational efficiency, it can be extended to composites containing a large number of particles. Particularly, a functionally graded material (FGM) can be treated as a bimaterial containing inhomogeneities with switched material phases and a continuously changing volume fraction  \cite{yin2004micromechanics,yin2007micromechanics,Wu2024_IJES}. The DR-iBEM provides a powerful tool for virtual experiments of FGMs. 

\subsection{Configuration of the FGM and thermal loading conditions}
Consider a ceramic-metal FGM sample, which aims to utilize the ceramic part for high thermal resistance and the metal part for high fracture toughness while avoiding significant thermal stresses at the bimaterial interface by a continuous material gradation \cite{Suresh1997}. The alumina (Al$_{\text{2}}$O$_\text{3}$) - nickel (NI) FGM composite will be studied with a linear gradation from Al$_{\text{2}}$O$_\text{3}$ at top surface to Nickel at bottom surface, whose dimensions are defined in Fig. \ref{fig:geometryc} that $l_a = l_b = h_1 = h_2 =$ 0.005 m. Assume that the initial temperature of the composite is uniform at $300$ K, and the boundary conditions are: (a) the temperature of the bottom surface is $300$ K; (b) the temperature of the top surface is set at $400$ K at $t>0$; and (c) all other surfaces are adiabatic. The transient heat transfer will gradually reach a steady state when $t \rightarrow \infty$. The thermal properties are as follows: (i) the thermal conductivity of Nickel and Al$_{\text{2}}$O$_\text{3}$ are $K^n = 90.7$ W / K $\cdotp$ m$^2$ and $K^a = 30.1$ W / K $\cdotp$ m$^2$ \cite{Suresh1997}, where superscripts $n$ and $a$ denote Nickel and Al$_{\text{2}}$O$_\text{3}$, respectively; (ii) the heat capacity of Nickel and Al$_{\text{2}}$O$_\text{3}$ are $C_p^n = 4.32 \times 10^6 \text{J} / \text{m}^3$ and $C_p^a = 3.96 \times 10^{6} \text{J} / \text{m}^3$ \cite{Chase1998} around the temperature 400 K. Symmetric properties of the FGM sample and boundary/initial conditions enables the simulation of a quarter of the sample, i.e., $x_1, x_2 \in [-0.0025, 0], x_3 \in [-0.005, 0.005]$ m.  

\begin{figure}
    \centering
    \includegraphics[width=0.4\linewidth]{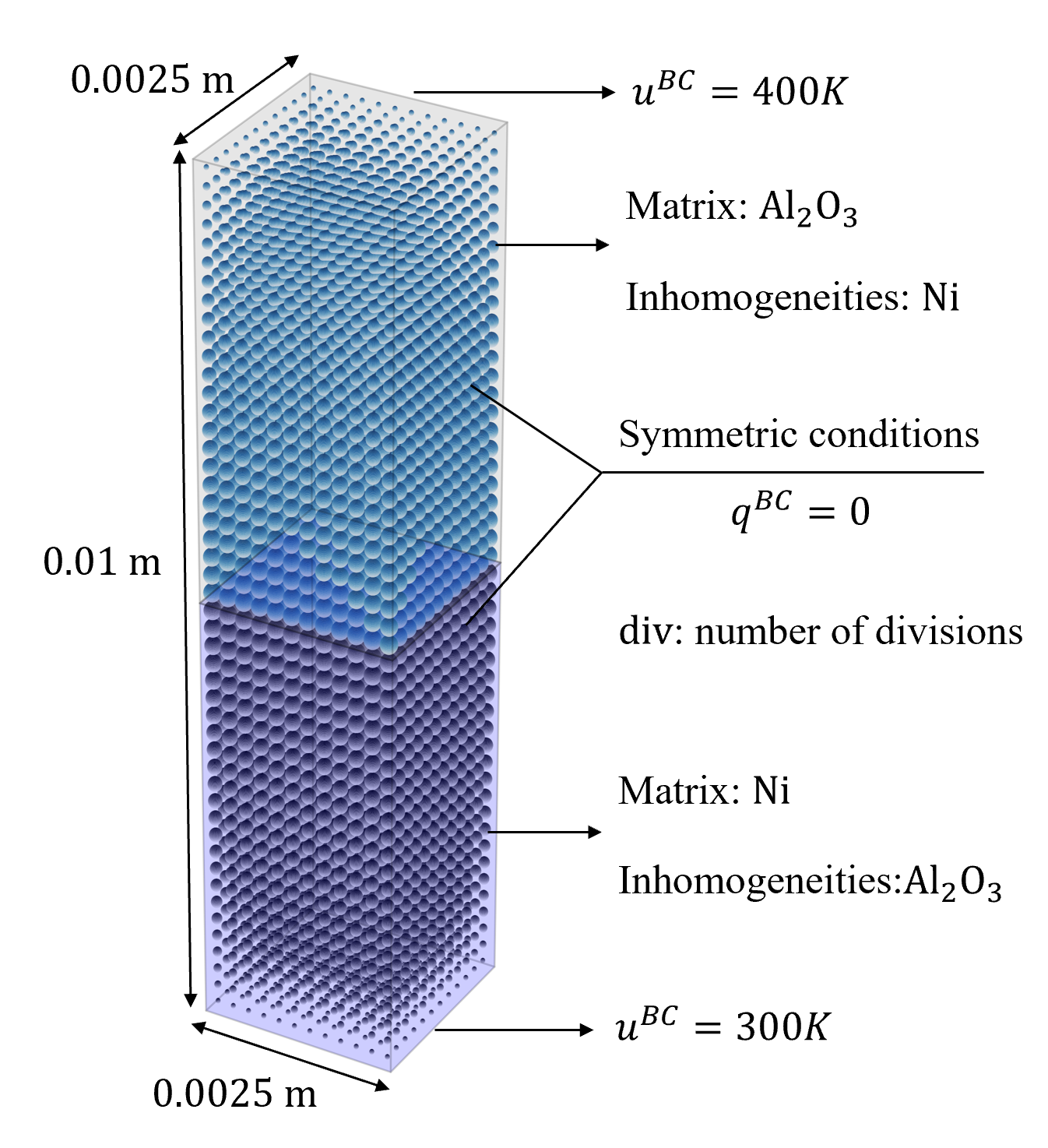}
    \caption{Schematic illustration of geometric and boundary conditions of a quarter of the FGM sample. A bimaterial model is utilized that: (i) in the upper domain ($\mathcal{D}^+$), NI particles filled in the Al$_{\text{2}}$O$_\text{3}$ matrix; (ii) in the lower domain ($\mathcal{D}^-$), Al$_{\text{2}}$O$_\text{3}$ particles filled in the Nickel matrix. ``div" refers to number of layers in the vertical direction to simulate gradation effects. Here, div = 48 and $6,912$ inhomogeneities are illustrated. }
    \label{fig:fgm_sample}
\end{figure}

To expedite the solving process and reduce computational cost, Fig. \ref{fig:fgm_sample} shows the bimaterial model \cite{Wu2024_IJES} that in the upper domain ($\mathcal{D}^+$), Al$_{\text{2}}$O$_\text{3}$ and NI serve as matrix and particle phases, respectively; whereas in the lower domain ($\mathcal{D}^-$), the matrix and particle phases are switched. Specifically, when the FGM exhibits linear gradation, the total volume fraction of inhomogeneities can be reduced from $50 \%$ to $25\%$. The gradation effects are simulated by finite layers (div) of spherical inhomogeneities, where the entire domain is uniformly divided into $\frac{\text{div}}{2}$ $\times$ $\frac{\text{div}}{2}$ $\times$ div elementary cubes, whose length is $\frac{2 l_a}{\text{div}}$. The center of each inhomogeneity is located at the center of each elementary cube, and the volume fraction change can be simulated by properly adjusting the radii of inhomogeneities along the gradation direction. For instance, Fig. \ref{fig:fgm_sample} shows the case div = 48, and $6,912$ inhomogeneities are involved and the maximum radius if $\frac{1}{4800} (3 / \pi)^{1/3} \approx 0.0001026$ m. The DR-iBEM is particularly suitable to predict the average thermal properties along the gradation direction and local field solutions \cite{Yin_iBEM}. The simulation time duration is $t \in [0, 20]$ s with the time interval $0.05$ s, and $1,802$ boundary nodes and $1,800$ quadrilateral boundary elements are utilized.

\subsection{Averaged temperature and heat flux along the gradation direction}
To evaluate the effective thermal properties of the FGM sample, the length ($x_1$ direction), width ($x_2$ direction), and height ($x_3$ direction) are uniformly divided into $10, 10, 20$ units, which correspondingly generates $10 \times 10 \times 20$ elementary cubes. Hence, $2,000$ sampling points are located at centers of cubes in one layer. For instance, when div = 48, totally $96,000$ sampling points are employed. Four cases of divisions, div = 20, 24, 32, 48 are studied with 500, 864, 2,048, and 6,912 inhomogeneities, respectively. Although the targeted volume fraction is linearly distributed along the gradation direction, the overall volume fraction of inhomogeneities decreases with the number of divisions due to the above discretization mechanism. Specifically, the volume fraction of particles in the upper phase changes among $0.325, 0.313, 0.297, 0.281$ accordingly \cite{Wu2024_IJES}. Therefore, the linear gradation can be better approximated by a larger number of divisions. 

\begin{figure}
\begin{subfigure}{.5\textwidth}
\centering
\includegraphics[width = 1\linewidth,height = \textheight,keepaspectratio]{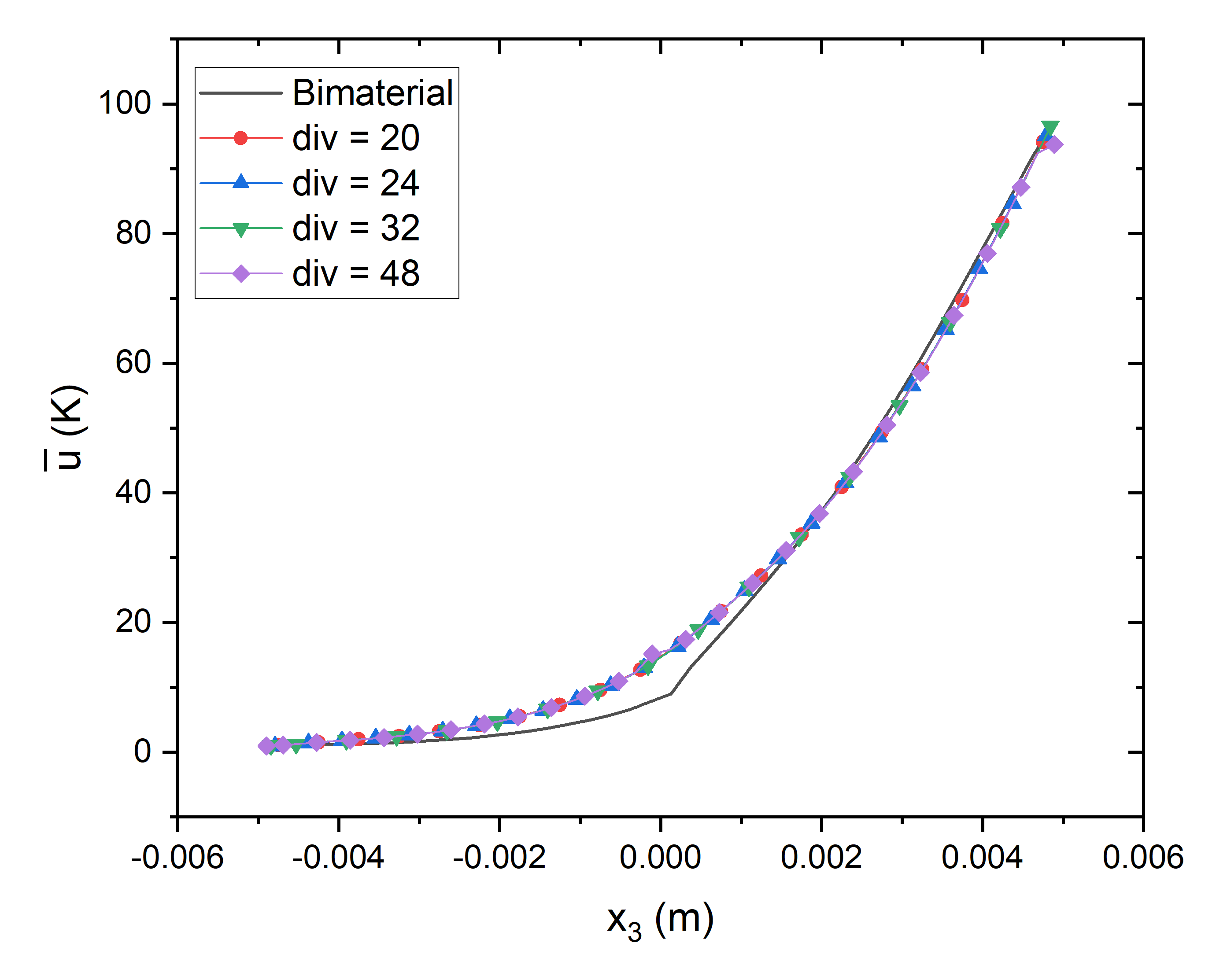}
\caption{}
\end{subfigure}
~
\begin{subfigure}{.5\textwidth}
\centering
\includegraphics[width = 1\linewidth,height = \textheight,keepaspectratio]{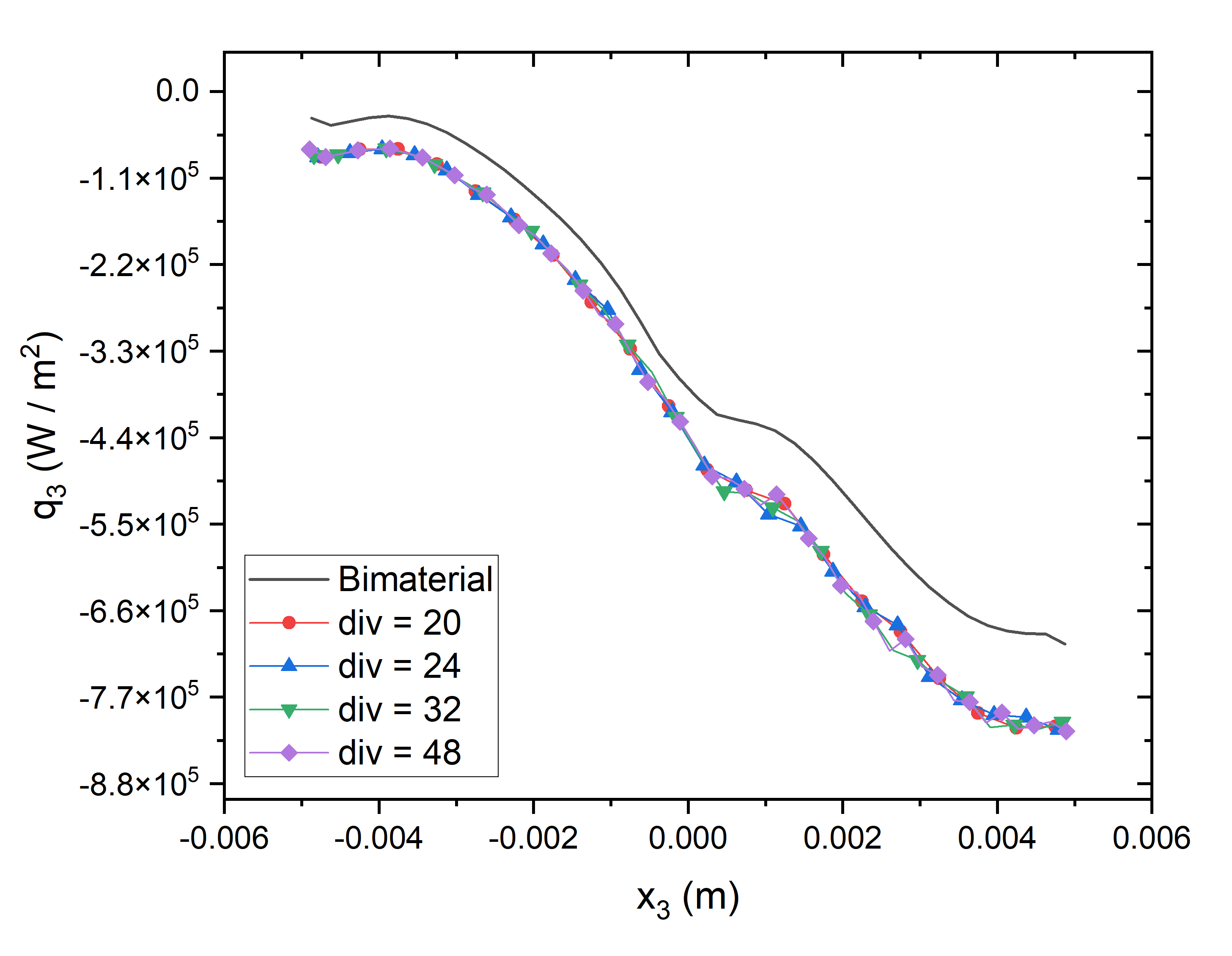}
\caption{}
\end{subfigure}
\caption{Variation and comparison of effective heat transfer behavior: (a) average temperature change $\tilde{u}$ and (b) heat flux $q_3$ at $0.8$ s of the FGM sample along the gradation direction for div = 20, 24, 32, and 48.}
\label{fig:AVG_0.8s}
\end{figure}

Fig. \ref{fig:AVG_0.8s} (a) plots the variations of averaged temperature (change) along the gradation direction at time 0.8s for different division cases. The legends ``div-20'', ``div-24'', ``div-32'', ``div - 48'' refer to the cases of div = 20, 24, 32, 48, respectively. Moreover, the legend ``Bimaterial'' represents the case of a bi-layered composite without inhomogeneities. Minor differences can be observed among the four cases. This observation agrees with our recent work \cite{Wu2024_IJES} (Fig. 17a,) where the flux boundary conditions were applied. Specifically, these four temperature curves are smooth across each layer, while the curve "Bimaterial" exhibits an obvious slope change at the bimaterial interface due to material mismatch. Due to the continuity conditions of the heat flux, Fig. \ref{fig:AVG_0.8s} (b) demonstrates that all curves are continuous throughout the plane $x_3 = 0$. Because the heat transfer process is unsteady at 0.8 s, while the heat gradually flows to points at lower heights, the maximum heat flux occurs at the top surface, and the minimum heat flux exists at the bottom surface. With the increase of time, the heat flux will approach a uniform distribution at the steady state. The component Ni dominates the lower phase, exhibiting higher thermal conductivity, which implies that temperature gradients are smaller at lower heights. Moreover, the curves associated with FGM composites are not as smooth as those in the bimaterial case, exhibiting fluctuations in each layer due to the disturbance caused by the inhomogeneities. Because the gradation is simulated by arrangements of multilayers of inhomogeneities, although the heat flux $q_3$ is supposed to be continuous along the gradation direction, there is no constraint on the slopes of heat flux variations. For instance, Fig. \ref{fig:thermal_spheroids} shows continuous distributions of heat flux $q_3$, but the local heat flux exhibits significant spatial variations within two inhomogeneities. As a result, it is rational for averaged thermal fields to show fluctuations due to the above approximated linear gradations. 

\begin{figure}
\begin{subfigure}{.5\textwidth}
\centering
\includegraphics[width = 1\linewidth,height = \textheight,keepaspectratio]{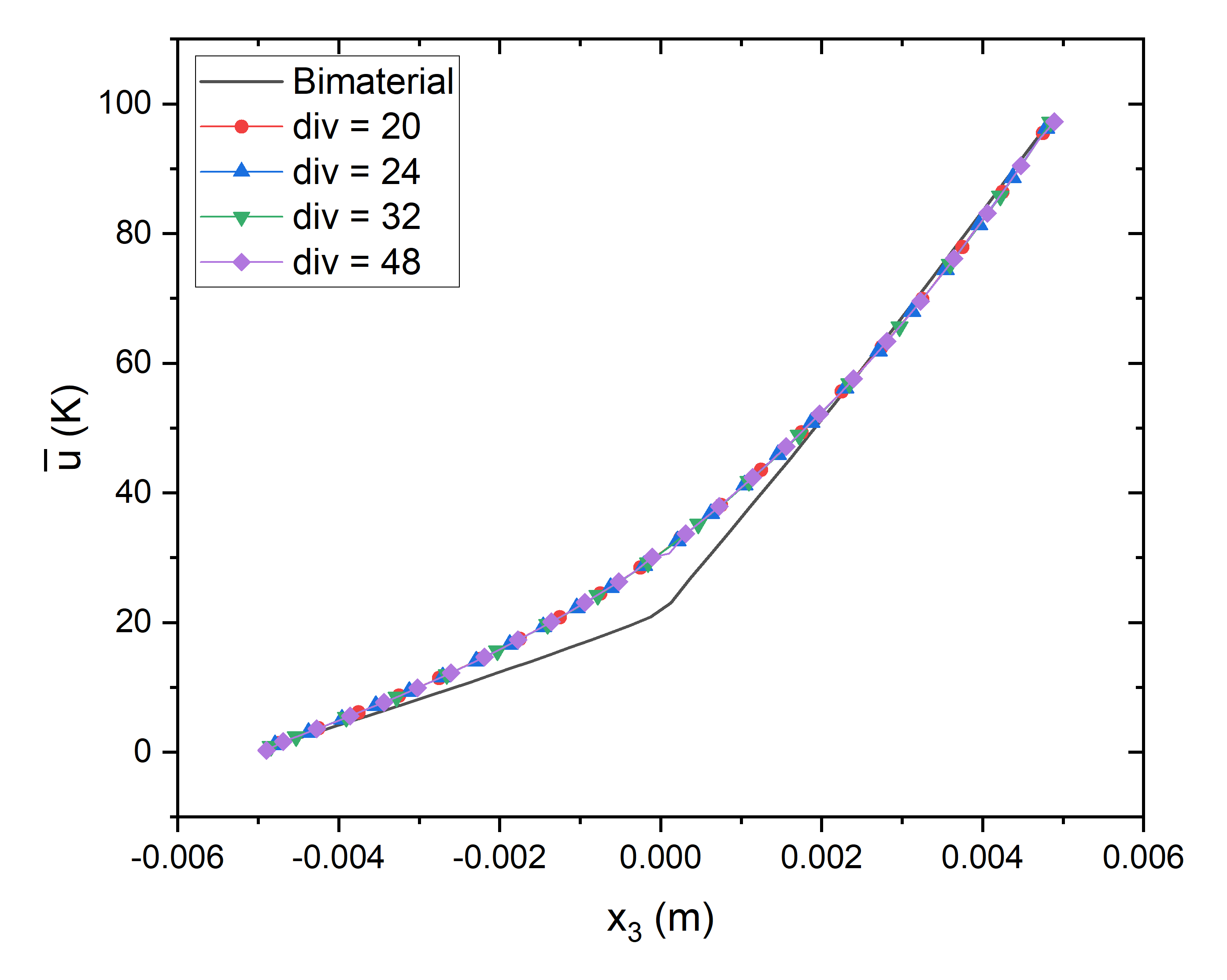}
\caption{}
\end{subfigure}
~
\begin{subfigure}{.5\textwidth}
\centering
\includegraphics[width = 1\linewidth,height = \textheight,keepaspectratio]{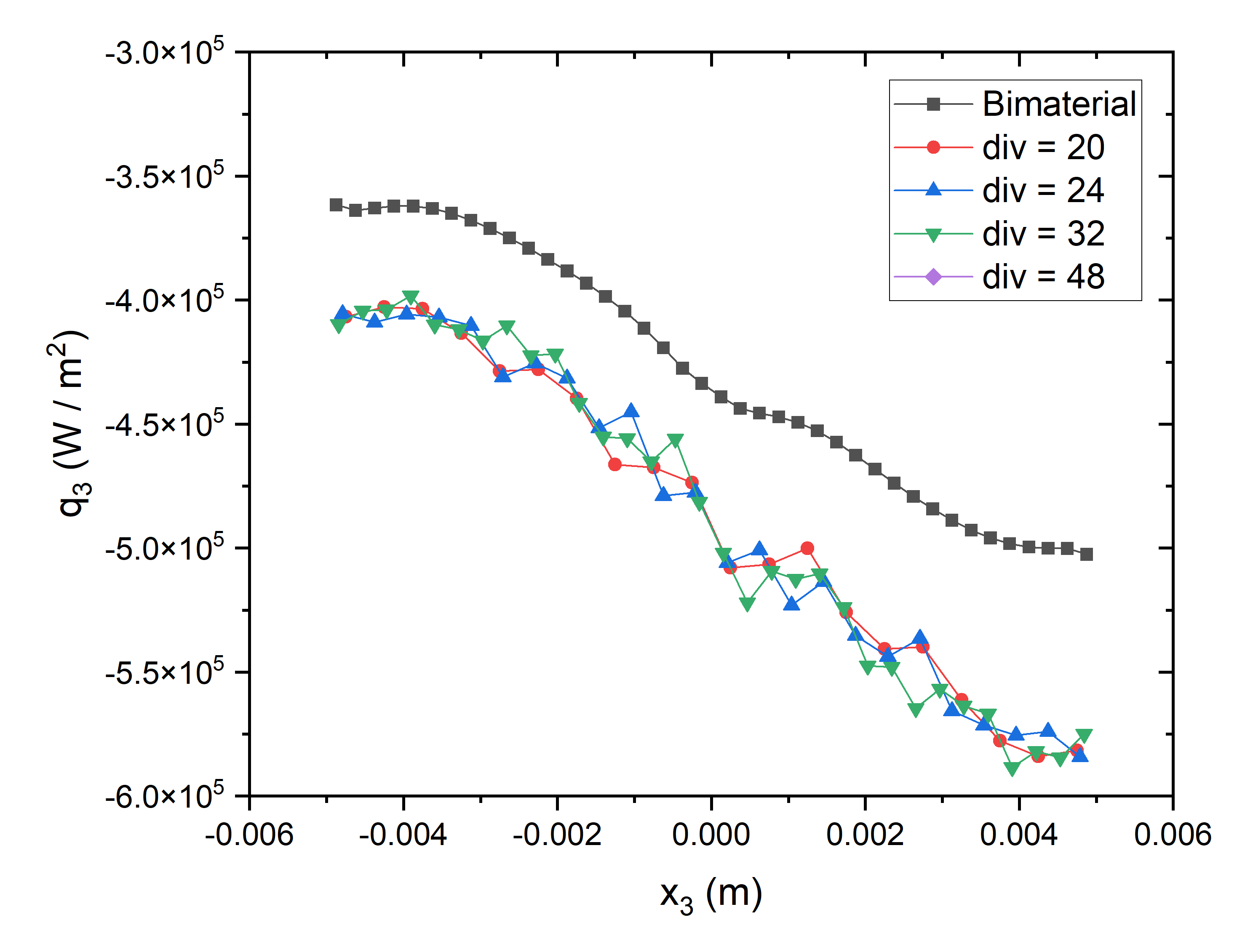}
\caption{}
\end{subfigure}
\caption{Variation and comparison of effective heat transfer behavior: (a) average temperature change $\tilde{u}$ and (b) heat flux $q_3$ at $1.5$ s of the FGM sample along the gradation direction for div = 20, 24, 32, and 48.}
\label{fig:AVG_1.5s}
\end{figure}

Subsequently, Figs. \ref{fig:AVG_1.5s} (a) and (b) show variations of the average temperature change and heat flux along the gradation direction at $1.5$s. Comparing curves in Fig. \ref{fig:AVG_1.5s} (a) with those in Fig. \ref{fig:AVG_0.8s} (a) shows that the average temperature changes at 1.5 s are greater as less heat is required to be stored in the material. For instance, the averaged temperature change at the bimaterial interface $x_3 = 0$ increases from $8.6$ to $21.5$ K, which implies that more heat has been transferred from the top surface to the bottom surface given more time. Similarly, the average temperature variations for FGM composites with different numbers of divisions are close to each other. Comparing Fig. \ref{fig:AVG_1.5s} (b) with  Fig. \ref{fig:AVG_0.8s} (b), we can see: (i) the range of heat flux becomes narrower, i.e., $q_3 \in [0.97, -8.42] \times 10^{5}$ $W / m^2$ in Fig. \ref{fig:AVG_0.8s} (b) and $q_3 \in [-4.13, -5.85] \times 10^{5}$ $W / m^2$ in Fig. \ref{fig:AVG_1.5s} (b), and the fluctuations among multiple layers become more obvious. When the ranges of heat flux become narrow,  the heat transfer process gradually shifts from a non-steady-state to a steady-state status.

\begin{figure}
\begin{subfigure}{.5\textwidth}
\centering
\includegraphics[width = 1\linewidth,height = \textheight,keepaspectratio]{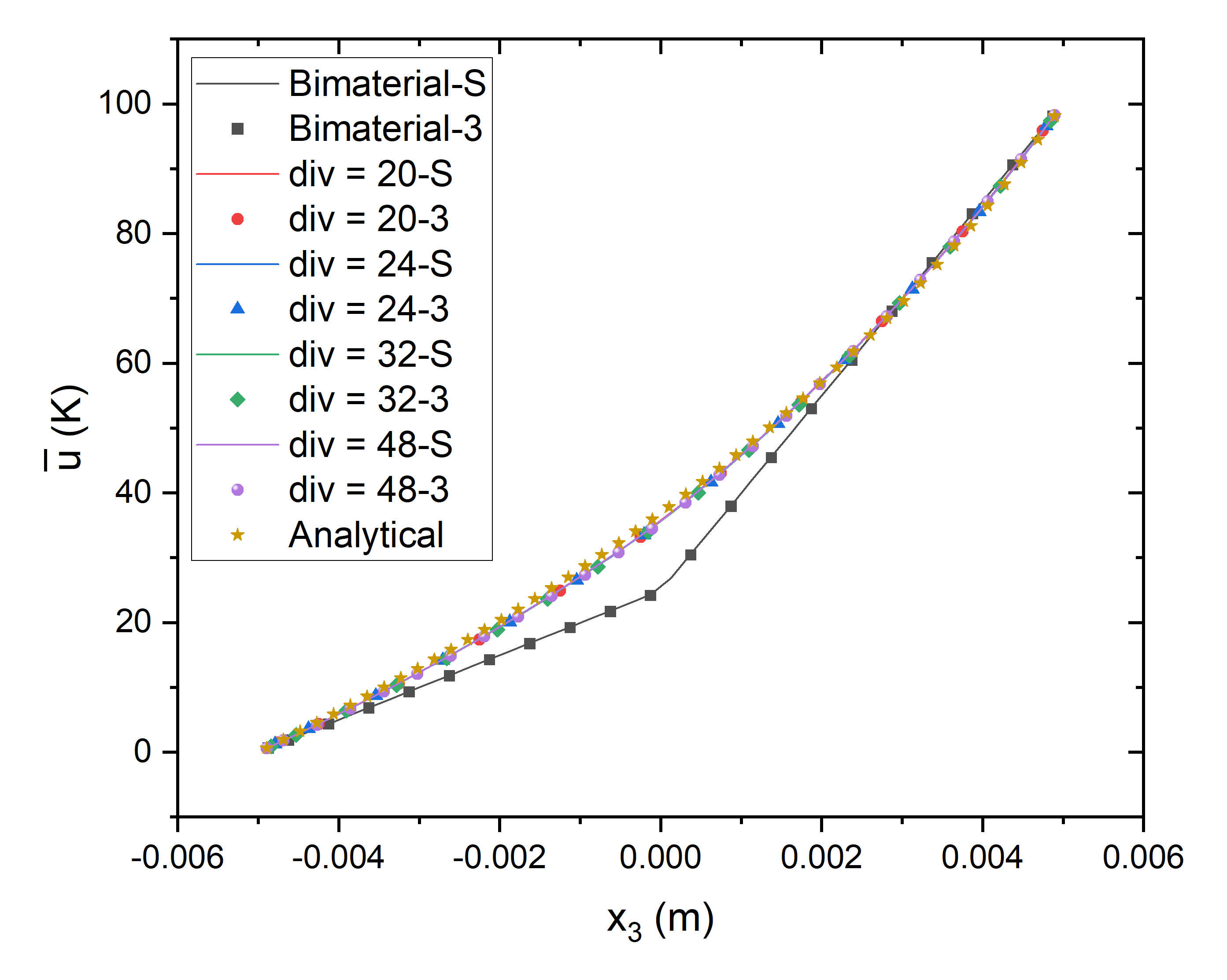}
\caption{}
\end{subfigure}
~
\begin{subfigure}{.5\textwidth}
\centering
\includegraphics[width = 1\linewidth,height = \textheight,keepaspectratio]{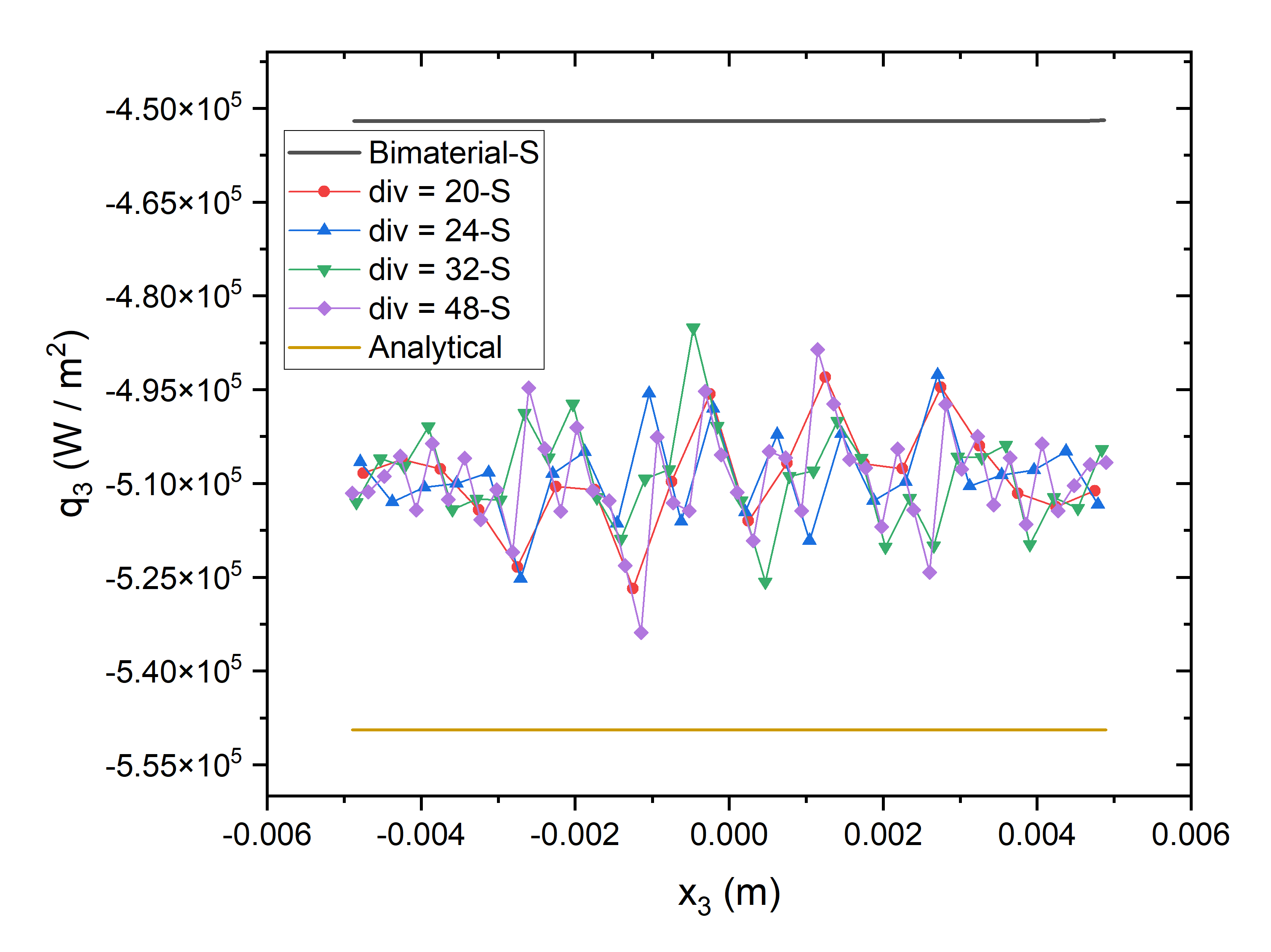}
\caption{}
\end{subfigure}
\caption{Variation and comparison of effective heat transfer behavior: (a) average temperature change $\tilde{u}$ and (b) heat flux $q_3$ at $3$ s of the FGM sample along the gradation direction for div = 20, 24, 32, and 48. Four numbers of division are involved, div = 20, 24, 32, 48.}
\label{fig:AVG_3s}
\end{figure}

When the time approaches a sufficiently long duration, i.e., 20 s, the transient heat transfer becomes steady-state heat conduction, and the results at 20 s match our recent work \cite{Wu2024_IJES} very well. Note that when the domain integral of the heat generation/storage rate becomes zero, the DR-iBEM for transient heat transfer indeed downgrades to iBEM for steady-state heat conduction. Fig. \ref{fig:AVG_3s} (a) compares variation of temperature (change) at 3 s with steady-state results, where curves ending with (3) and (S) refer to results at 3 s and steady-state condition, respectively. In addition, Fig. \ref{fig:AVG_3s} (a) shows that temperature curves at $3$ s exhibit minor discrepancies from the results at the steady state. For instance, the maximum difference is less than $0.4 \%$. Therefore, thermal fields at 3 s can be approximately considered as the steady-state status, and the temperature change over time can be very minor for $t>3$ s. 
 
Comparing the steady-state temperature fields for $t>3$ s with those in Figs. \ref{fig:AVG_0.8s} (a) and \ref{fig:AVG_1.5s} (a), we can observe:: (i) For the bimaterial case, the temperature distribution is linear and continuous with a sudden slope change at the interface, which is caused by discontinuous temperature gradients owing to dissimilar thermal conductivities. However,  when the heat transfer is not in the steady state, the temperature curves along the gradation direction are not linear, as observed in Figs. \ref{fig:AVG_0.8s} (a) and \ref{fig:AVG_1.5s} (a). (ii) Regarding the FGM cases, the thermal conductivity changes continuously along the gradation direction, which results in a variation of the slope in each layer and shows that the FGM composites exhibit smooth distributions of average temperature. Overall, the averaged temperature curves for FGMs are nonlinear under both unsteady and steady heat conduction cases. 

Fig. \ref{fig:AVG_3s} (b) plots the steady-state heat flux variations along the gradation direction. For the bimaterial case, the heat flux is constant, which leads to piece-wise temperature gradients in the steady state shown in Fig. \ref{fig:AVG_3s} (a). Regarding the FGM samples, the range of heat flux changes is narrow at the steady-state status, leading to greater local fluctuations than in the previous unsteady-state cases, as more heat transfer through the FGM. Essentially, the FGM sample can be homogenized into a one-dimensional composite along the gradation direction. Let $T' - T''$ refer to the temperature difference between the top and bottom surfaces. The analytical solutions for the steady-state heat transfer in the bimaterial case and FGM case can be derived by solving the governing equation $\nabla (K(\textbf{x}) \nabla T(x)) = 0$. Specifically, (i) for the bimaterial case, the steady-state heat flux is $\frac{-2 (T' - T'') K' K''}{(h_1 + h_2) (K' + K'')} \approx -4.5 \times 10^5 W / m^2$, which agrees well with the ``Bimaterial'' curve in Fig. \ref{fig:AVG_3s} (b); (ii) for FGM with linear gradation, the temperature distribution is $(T' - T'') \text{log} \frac{2K'L}{(3K'- K'')L + 2(K'' - K') x_3} / \text{log} \frac{K'}{K''}$, and the heat flux distribution is $(T' - T'')\frac{K' - K''}{(h_1 + h_2) \text{log} \frac{K''}{K'}}$. Fig. \ref{fig:AVG_3s} (a) shows that the temperature evaluated for the FGM with linearly distributed thermal properties is close to the DR-iBEM predictions with a deterministic microstructure, and the maximum difference exhibits around the plane $x_3 = 0$. Fig. \ref{fig:AVG_3s} (b) illustrates the difference between the analytical solution and DR-iBEM prediction. Note that the analytical solution assumes that the thermal conductivity linearly changes along the gradation direction, i.e., from K' to K'', $K(\textbf{x}) = K' + (K' - K'') \frac{x_3}{h_1 + h_2}$. However, the linear change of volume fractions is not equivalent to the linear gradation of thermal conductivity with equal divisions. Although such a solution has good predictions on temperature fields, it overestimates the heat flux. The averaged heat flux is approximately predicted by DR-iBEM at $-5.10 \times 10^{5} W / m^2$ is in between the results of the bimaterial case and analytical solution for the linear-graded FGM case, which is $0.13$ times greater than the average heat flux of the bimaterial case. The heat flux is highly sensitive to the distribution of thermal properties, so high-fidelity predictions of effective thermal properties are required. 

\subsection{Local temperature and heat flux of FGM samples}
The key difference between a homogenized material and the actual material sample with a certain microstructure is that the latter can provide accurate predictions of local fields from the solution of a boundary value problem. Many realizations of the microstructural solution can predict a statistical prediction for the homogenized material. 
To illustrate DR-iBEM's capabilities solving local fields, contour plots on temperature and heat flux are provided as follows: (i) the plane $x_1 = -0.00125$ m is selected, while $x_2$ and $x_3$ are variables; (ii) the lengths along $x_2, x_3$ directions are evenly divided into $160$ and $640$ parts, respectively, which forms $160 \times 640$ elementary cubes with $102,400$ sampling points located at the centers of elementary cubes. Note that the previous subsection demonstrates that when $t = 3$ s, the transient heat transfer process has almost shifted to steady-state heat conduction. Therefore, this subsection displays contour plots at time 0.8, 1.5, and 3 s, while div = 20 and 36. 

\begin{figure}
\begin{subfigure}{.3\textwidth}
\centering
\includegraphics[width = 1\linewidth,height = \textheight,keepaspectratio]{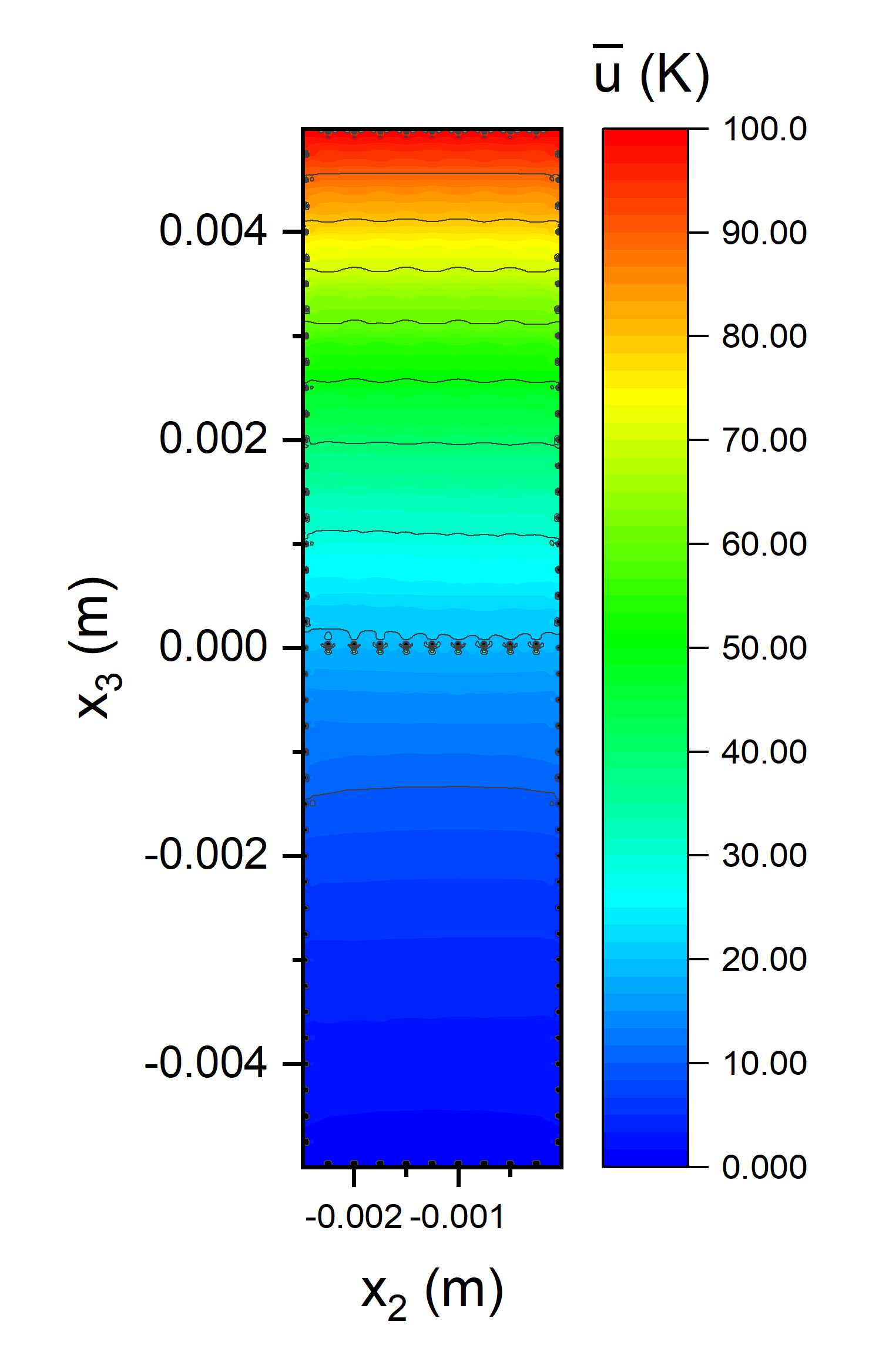}
\caption{}
\end{subfigure}
~
\begin{subfigure}{.3\textwidth}
\centering
\includegraphics[width = 1\linewidth,height = \textheight,keepaspectratio]{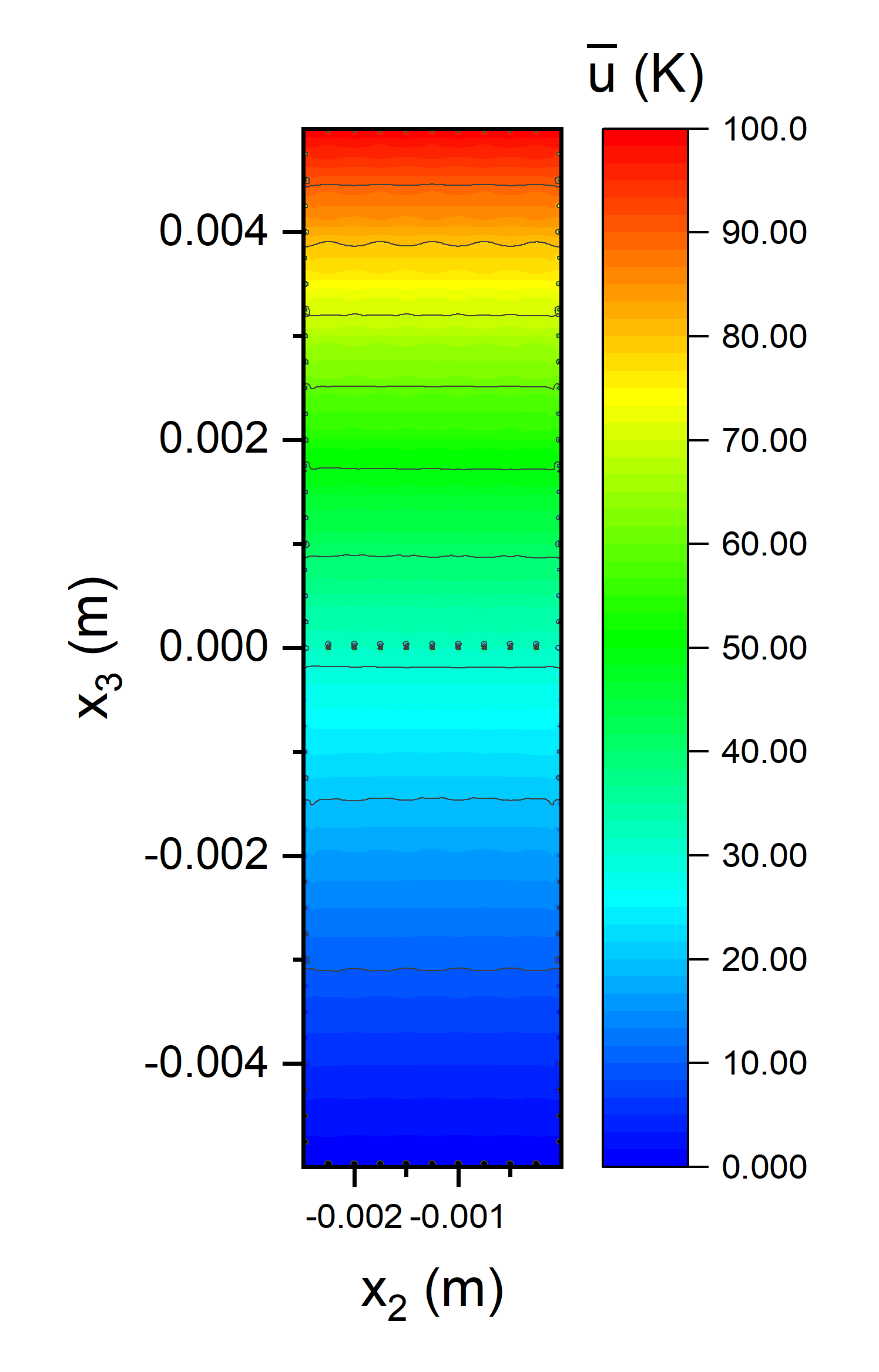}
\caption{}
\end{subfigure}
~
\begin{subfigure}{.3\textwidth}
\centering
\includegraphics[width = 1\linewidth,height = \textheight,keepaspectratio]{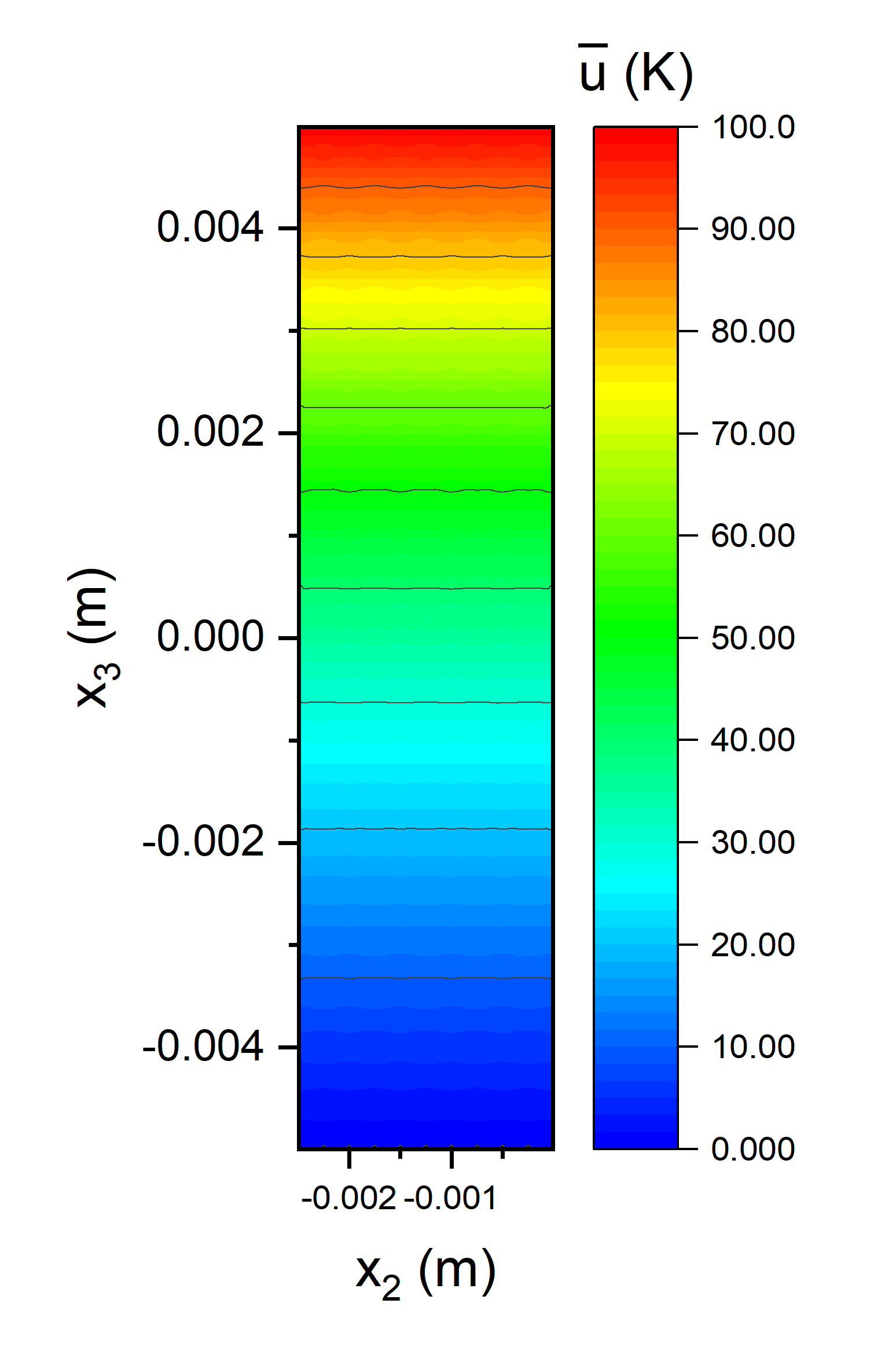}
\caption{}
\end{subfigure}

\begin{subfigure}{.3\textwidth}
\includegraphics[width = 1\linewidth,height = \textheight,keepaspectratio]{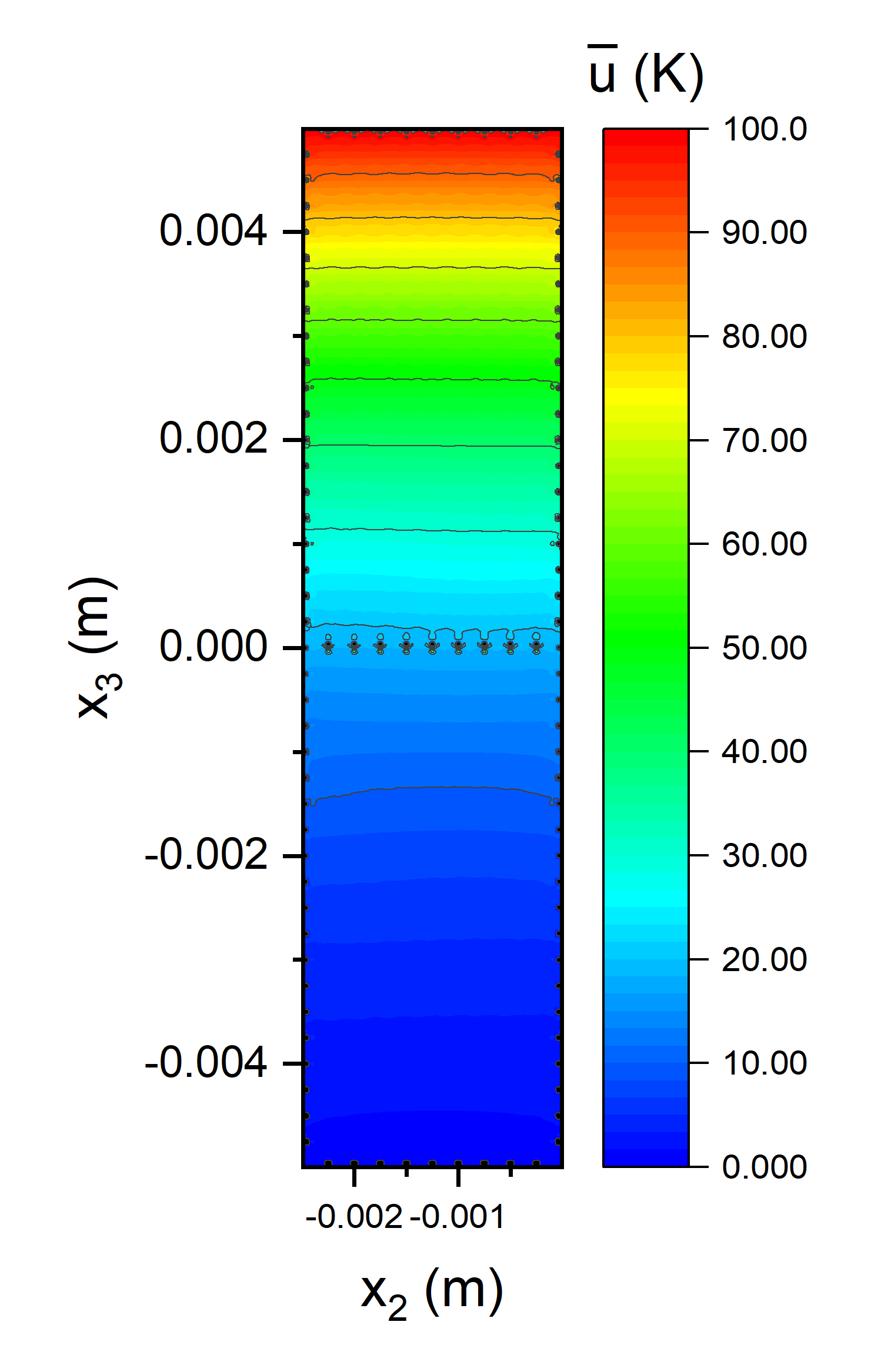}
\caption{}
\end{subfigure}
~
\begin{subfigure}{.3\textwidth}
\centering
\includegraphics[width = 1\linewidth,height = \textheight,keepaspectratio]{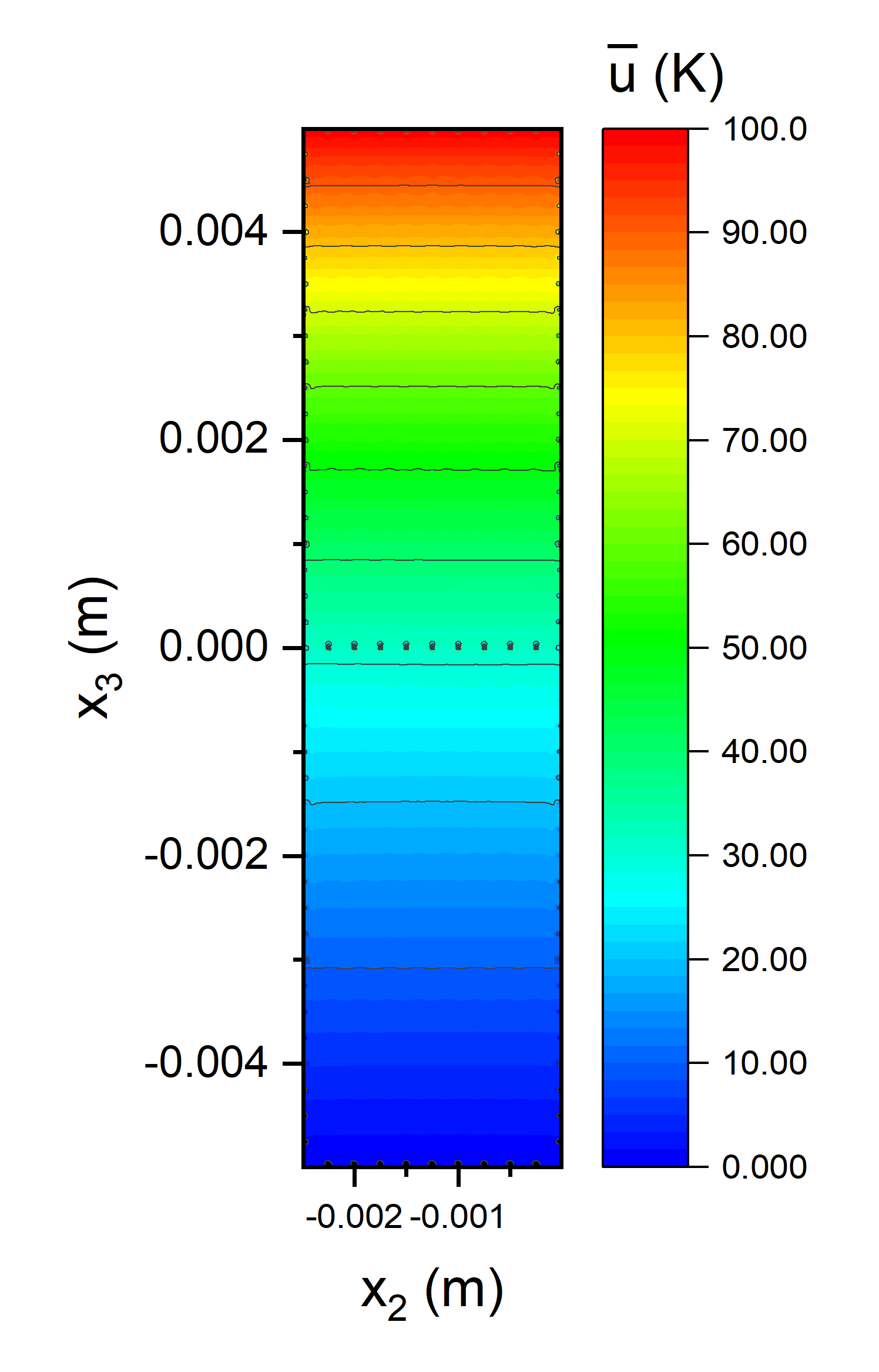}
\caption{}
\end{subfigure}
~
\begin{subfigure}{.3\textwidth}
\centering
\includegraphics[width = 1\linewidth,height = \textheight,keepaspectratio]{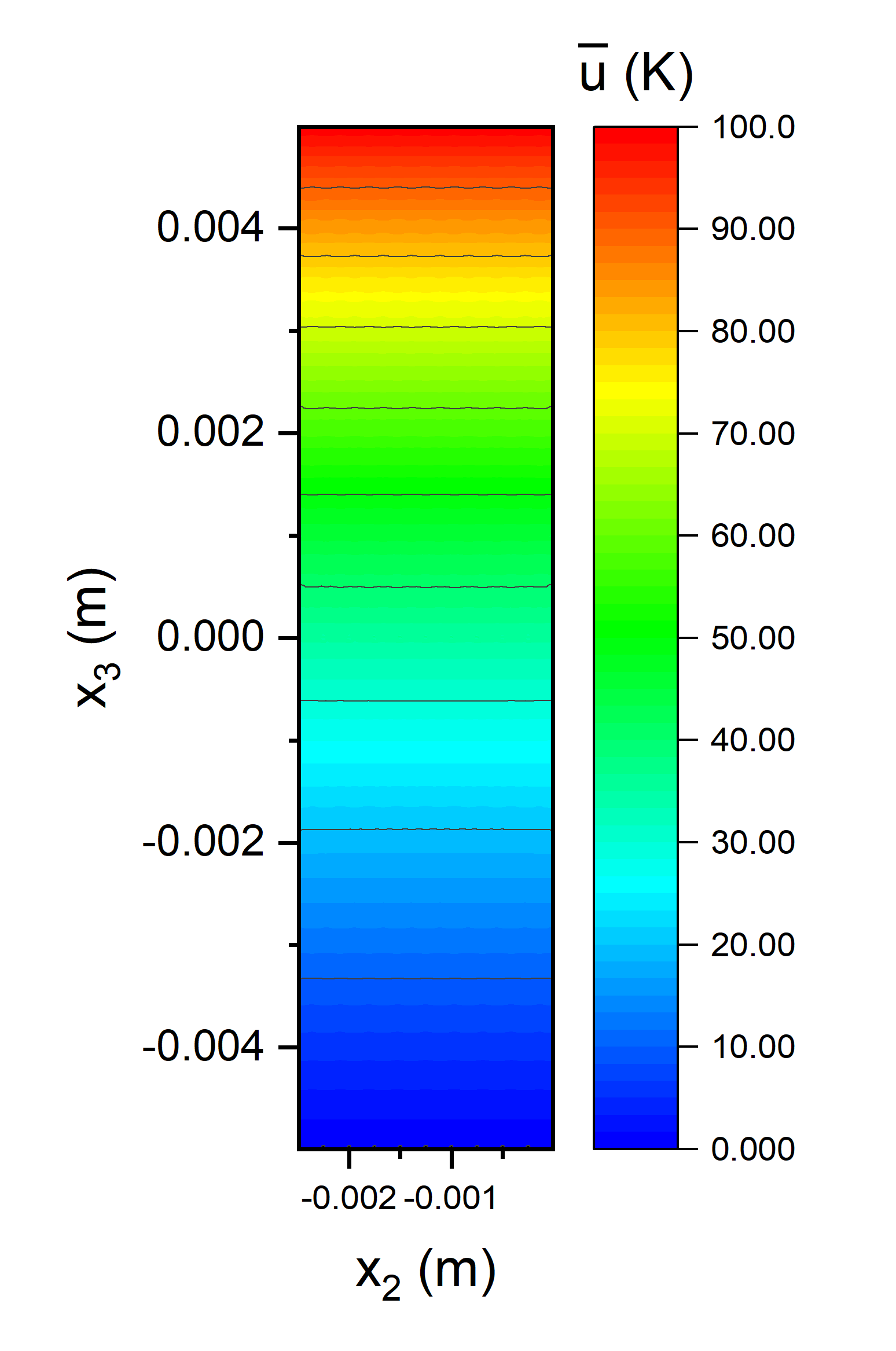}
\caption{}
\end{subfigure}

\caption{Temperature change contour of the internal cross-section of a quarter FGM sample ($x_2 \in [-0.0025, 0]$ m, $x_3 \in [-0.005, 0.005]$ m). Sub-figures (a), (b), (c) plot temperature contour when div = 20 at time 0.8, 1.5, 3s, respectively; Sub-figures (c), (d), (e) plot temperature contour when div = 36 at time 0.8, 1.5, 3s, respectively.}
\label{fig:contour_temp}
\end{figure}

Figs. \ref{fig:contour_temp} (a-c) plot the contour of temperature change when div = 20 at time 0.8, 1.5, 3s, respectively. Comparing Fig. \ref{fig:contour_temp} (a) to Fig. \ref{fig:contour_temp} (b) shows that the heat flow takes time to transfer from the top surface to the bottom surface. Therefore, Fig. \ref{fig:contour_temp} (b) exhibits higher temperature for all internal points, which has been shown in Figs. \ref{fig:AVG_0.8s} (a) and \ref{fig:AVG_1.5s} (a). Note that in Figs. \ref{fig:contour_temp} (a-c), it is difficult to observe shapes of inhomogeneities, although the contour plots show the temperature distribution of the cross-section. Such a phenomenon can be well explained using the steady-state Green's function, which suggests that the temperature should be continuous at all internal points. Therefore, temperature typically does not exhibit significant variations at entering the region of inhomogeneities. However, contours in Figs. \ref{fig:contour_temp} (a) and Fig. \ref{fig:contour_temp} (b) have small curves, which can be observed from the contour lines around $x_3 = 0.004$ m. This small curvature reveals disturbances caused by inhomogeneities, and such effects are more obvious for internal points closer to the top surfaces. Note that although equivalent conditions are consistent across all inhomogeneities, regardless of their location in either the upper or lower matrix phases, the disturbances are determined by the magnitude of eigen-fields. Based on Eq. (\ref{eq:equiv_flux_transient}) and Eq. (\ref{eq:equiv_heat_transient}), the magnitude of eigen-fields are highly related to unperturbed thermal fields. For unsteady-state heat transfer, the temperature change is primarily driven by heat flux from the upper surface, which induces greater unperturbed thermal fields, leading to larger eigenfields and more pronounced disturbances close to the upper surface. Moreover, the inhomogeneities in the neighborhood of the plane $x_3 = 0$ have larger radii, and the disturbances can be even higher. As the heat transfer process approaches the steady state, the contours exhibit smaller curvatures as the unperturbed heat flux becomes narrower. Consequently, although the disturbances from inhomogeneities still exist, inclusions serve as source fields with approaching values symmetric to the contour line, which flatten the contours. Keep in mind that similar patterns of temperature distribution can be found in Figs. \ref{fig:contour_temp} (d-f) when div = 36. When div is greater, inhomogeneities are comparatively smaller; therefore, contour lines in Figs. \ref{fig:contour_temp} (d-e) has more small curves. Although inhomogeneities exhibit certain size effects, the temperature distribution issimilar between div = 20 and div = 36, which agrees with Figs. \ref{fig:AVG_0.8s} (a), \ref{fig:AVG_1.5s} (a), and \ref{fig:AVG_3s} (a). 

\begin{figure}
\begin{subfigure}{.3\textwidth}
\centering
\includegraphics[width = 1\linewidth,height = \textheight,keepaspectratio]{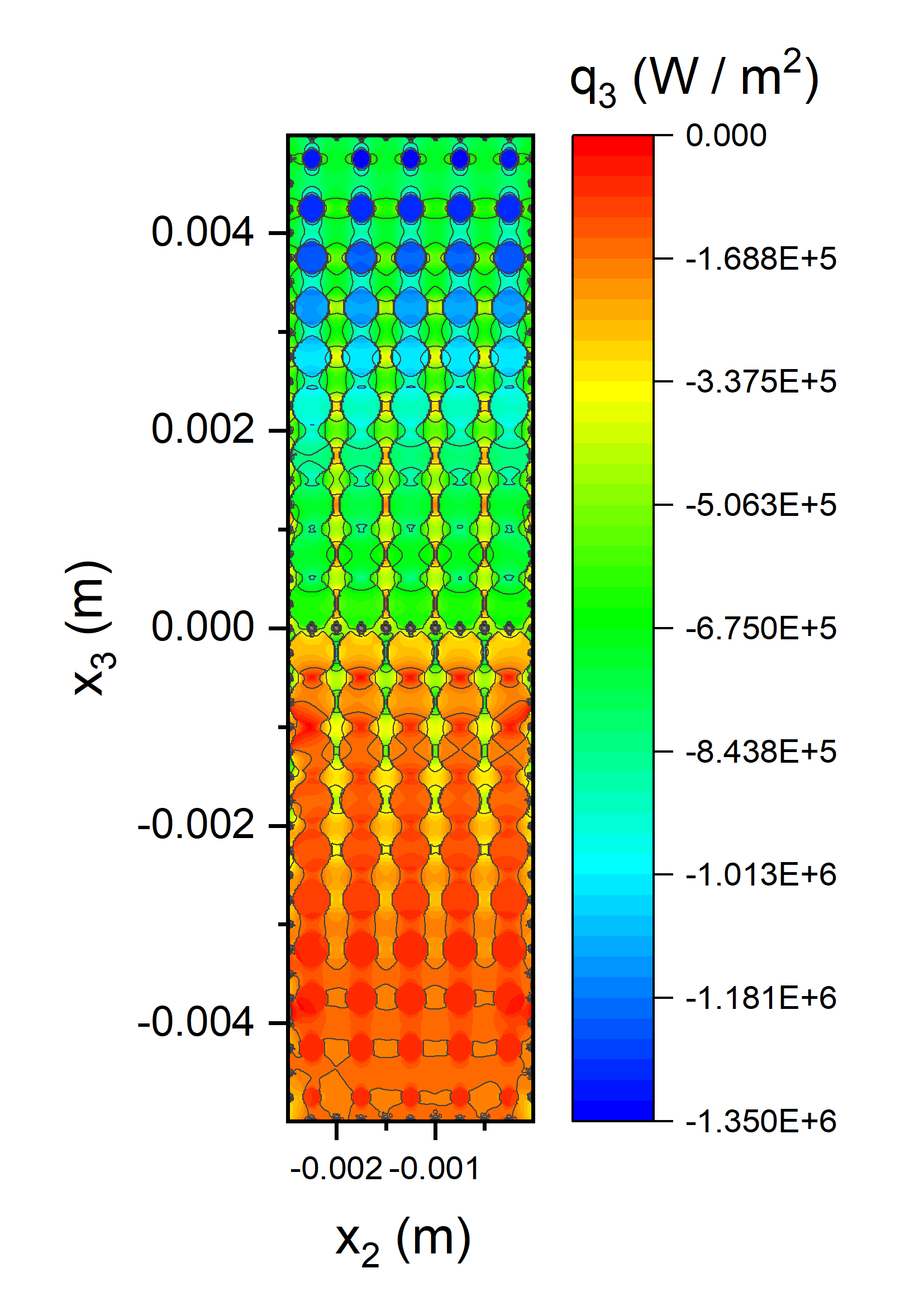}
\caption{}
\end{subfigure}
~
\begin{subfigure}{.3\textwidth}
\centering
\includegraphics[width = 1\linewidth,height = \textheight,keepaspectratio]{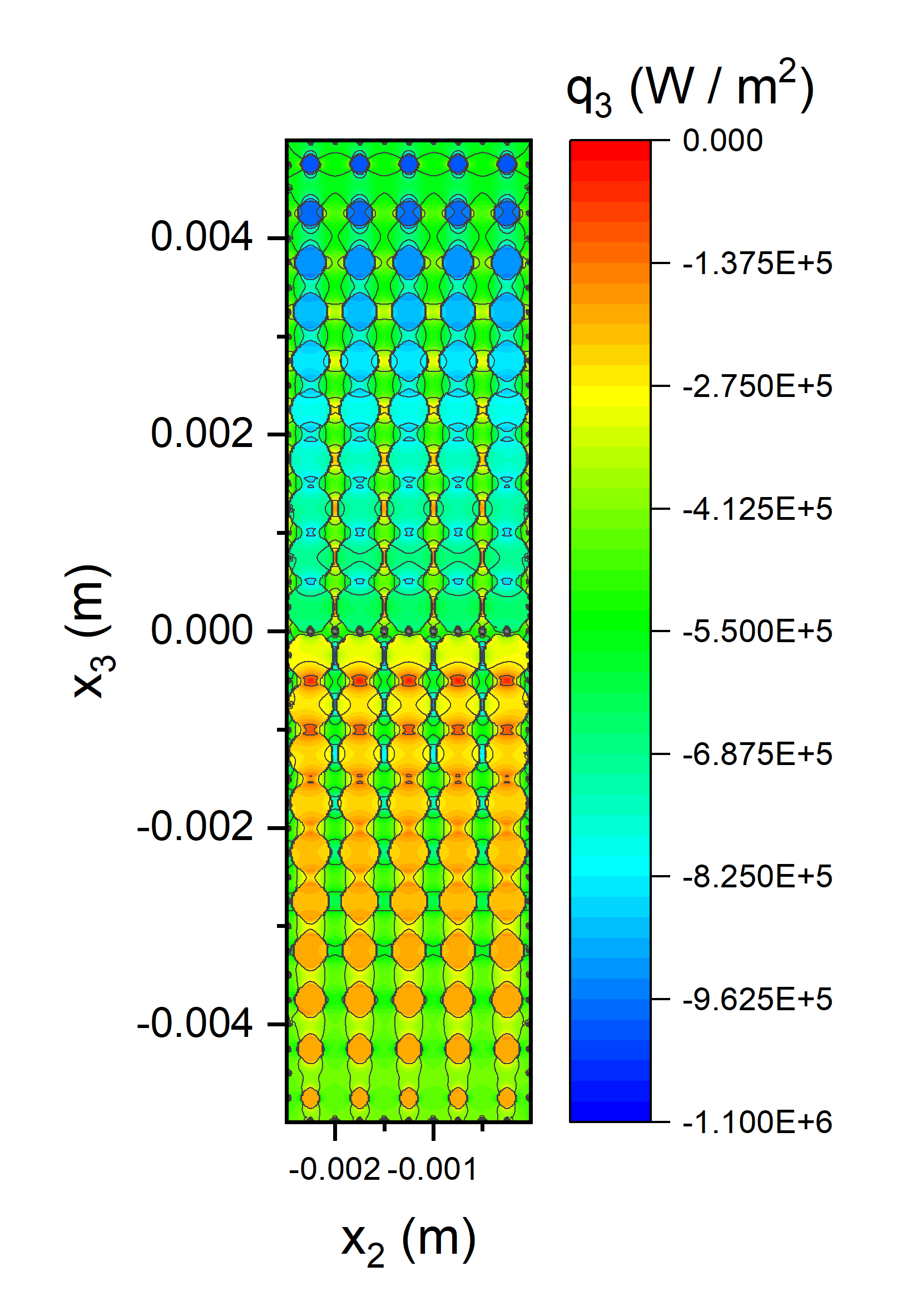}
\caption{}
\end{subfigure}
~
\begin{subfigure}{.3\textwidth}
\centering
\includegraphics[width = 1\linewidth,height = \textheight,keepaspectratio]{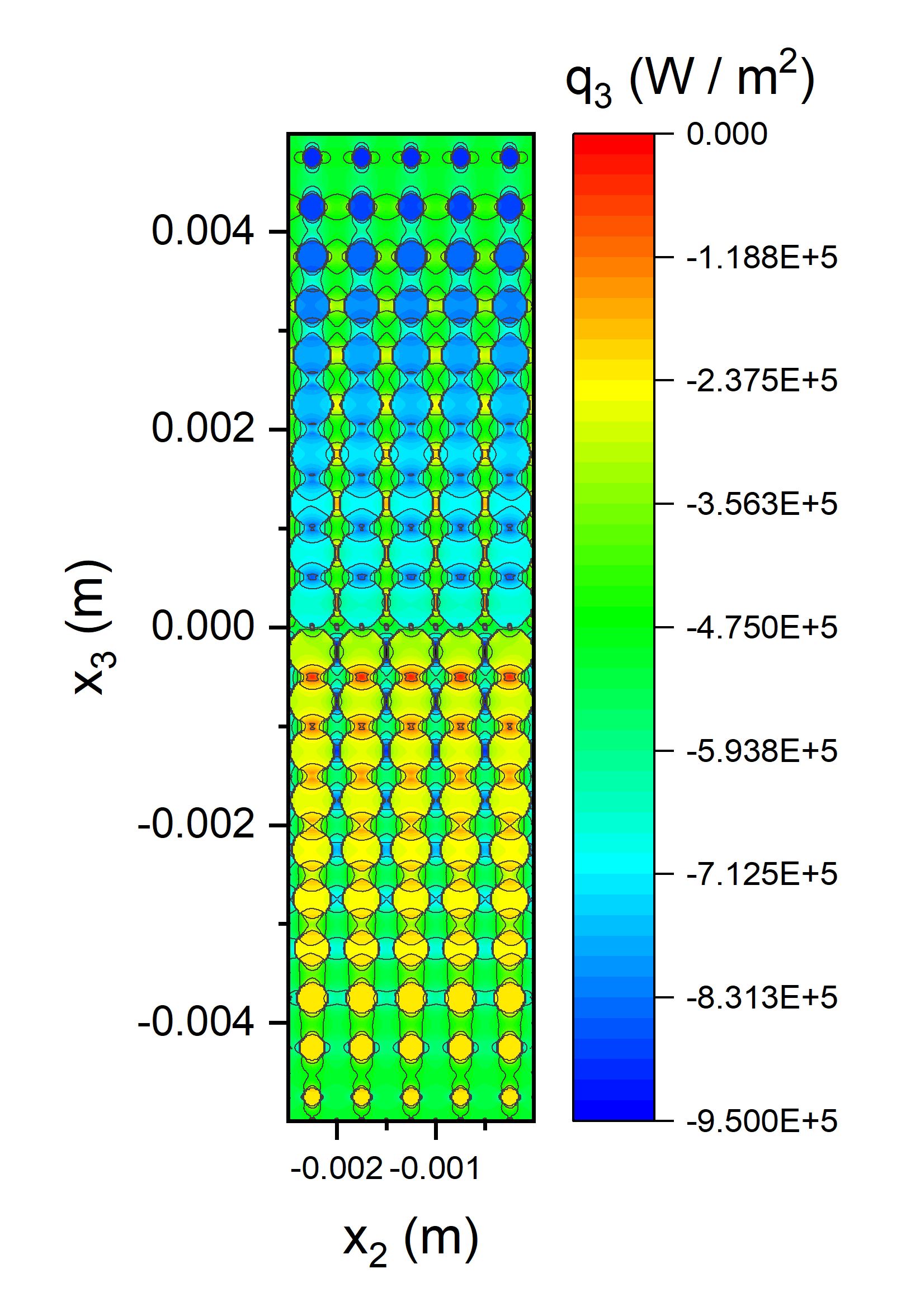}
\caption{}
\end{subfigure}

\begin{subfigure}{.3\textwidth}
\includegraphics[width = 1\linewidth,height = \textheight,keepaspectratio]{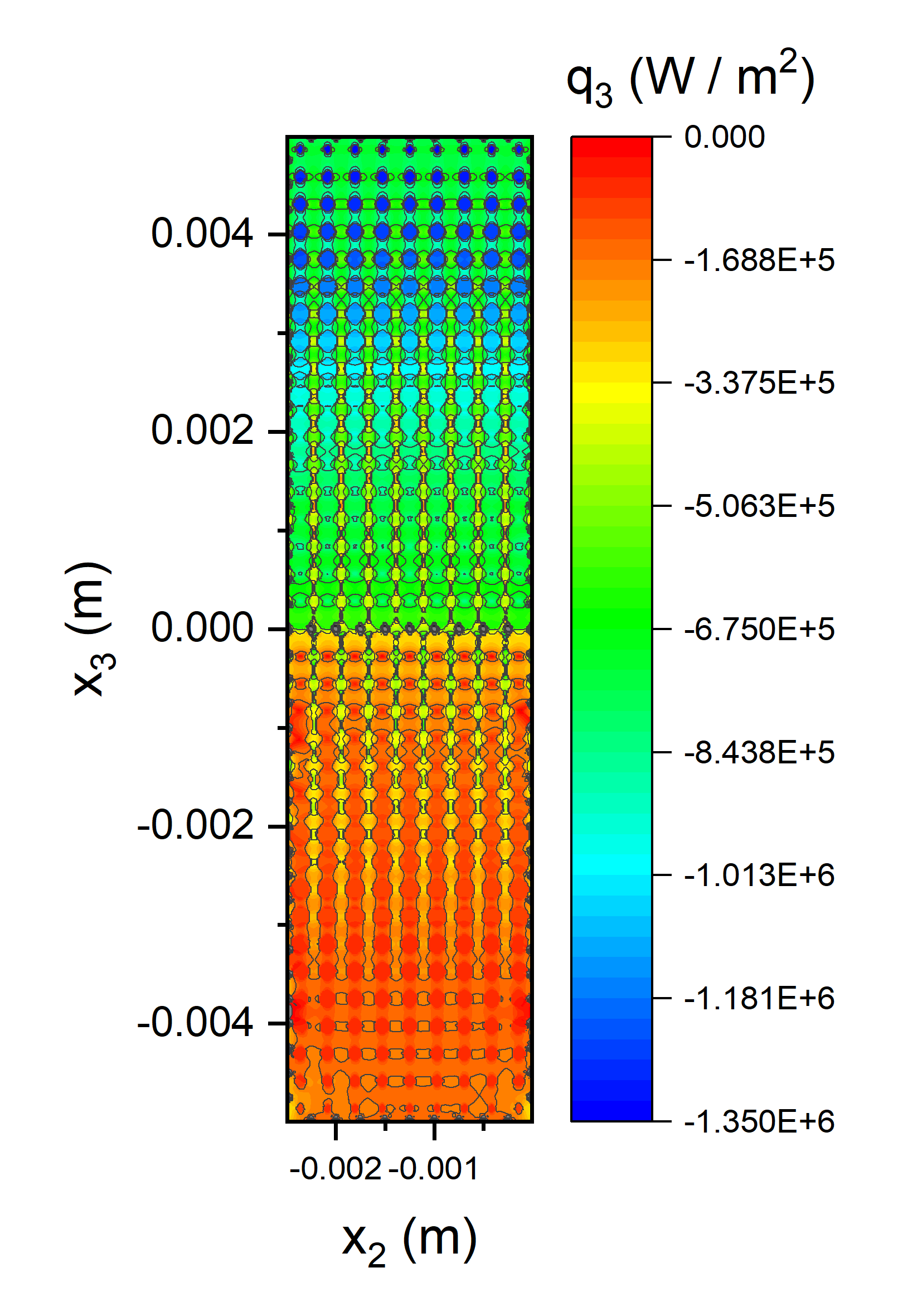}
\caption{}
\end{subfigure}
~
\begin{subfigure}{.3\textwidth}
\centering
\includegraphics[width = 1\linewidth,height = \textheight,keepaspectratio]{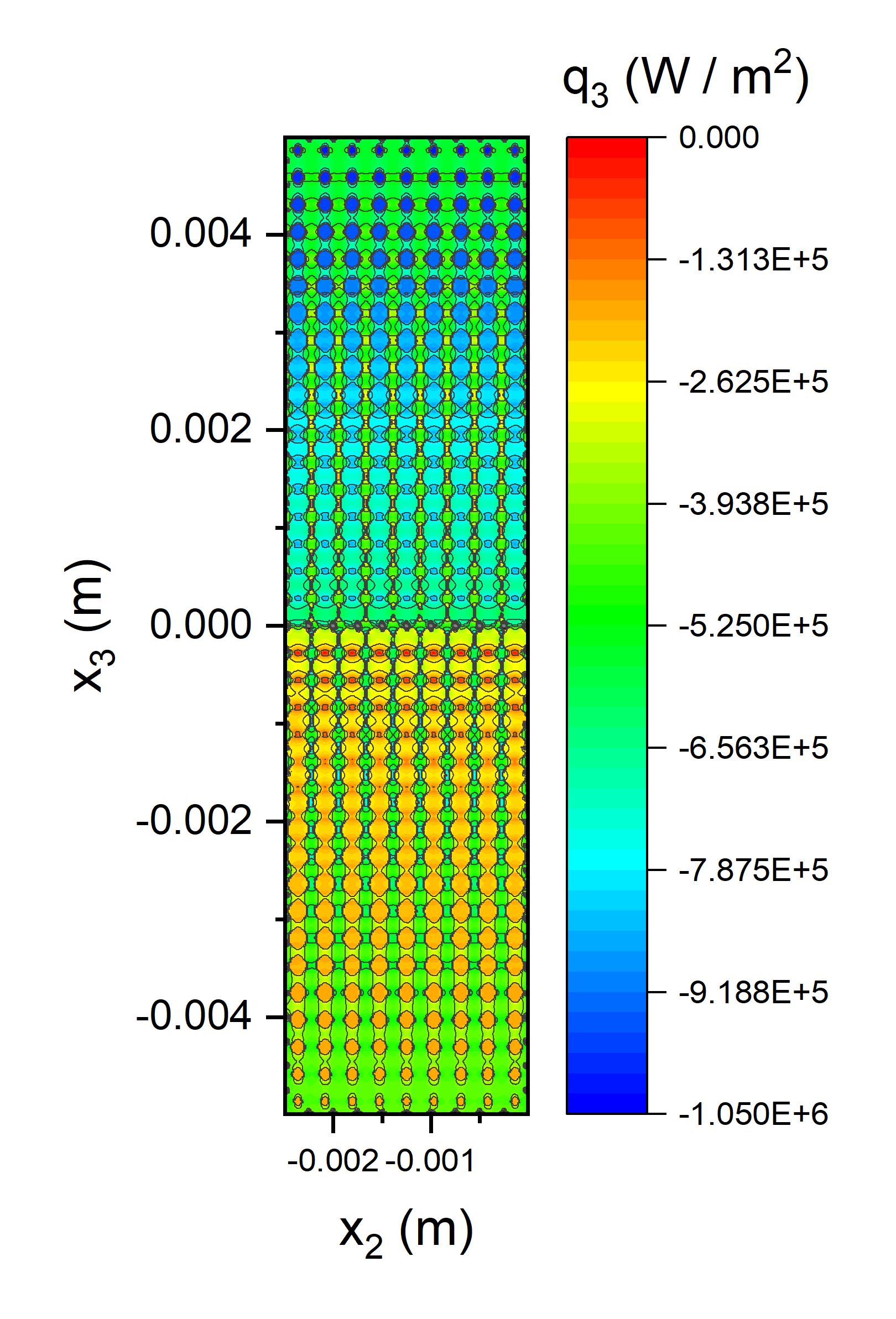}
\caption{}
\end{subfigure}
~
\begin{subfigure}{.3\textwidth}
\centering
\includegraphics[width = 1\linewidth,height = \textheight,keepaspectratio]{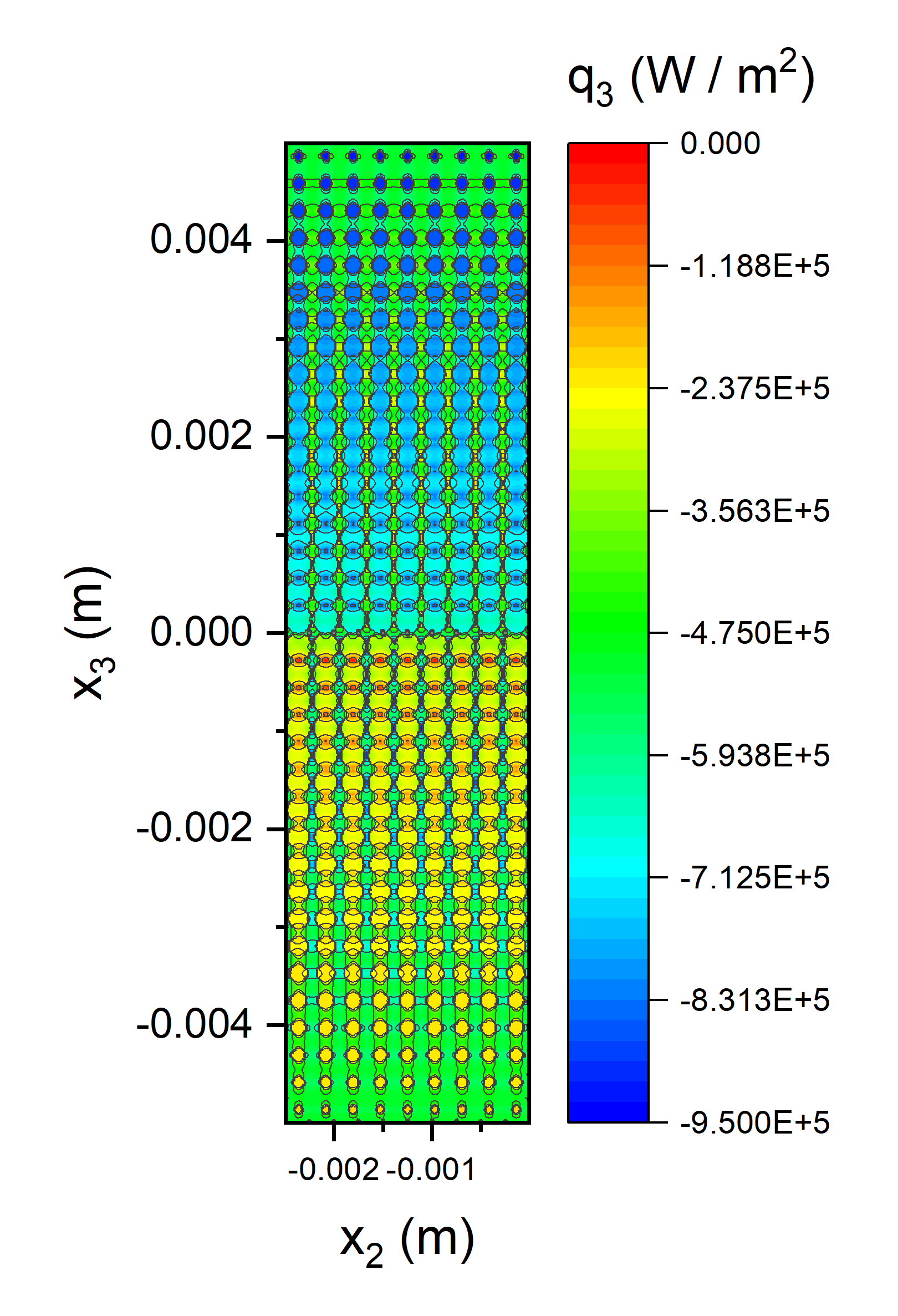}
\caption{}
\end{subfigure}

\caption{Heat flux contour of the internal cross-section of a quarter FGM sample ($x_2 \in [-0.0025, 0]$ m, $x_3 \in [-0.005, 0.005]$ m). Sub-figures (a), (b), (c) plot flux contour when div = 20 at time 0.8, 1.5, 3s, respectively; Sub-figures (c), (d), (e) plot flux contour when div = 36 at time 0.8, 1.5, 3s, respectively.}
\label{fig:contour_flux}
\end{figure}

Compared Figs. \ref{fig:contour_flux} (a-c), the ranges of heat flux become narrower, which was illustrated in Fig. \ref{fig:AVG_0.8s} (b) and Fig. \ref{fig:AVG_3s} (b). Unlike the previous temperature contours in Fig. \ref{fig:contour_temp}, it is straightforward to identify positions and size of inhomogeneities of the cross-section. Such a phenomenon can be explained by Eshelby's tensors in Eq. (\ref{eq:Eshelby_tensors}), which associates with domain integrals of two potential functions. Although the heat flux $q_3$ is continuous for surface normals (0, 0, -1) or (0, 0, 1), the heat flux within the inhomogeneity can exhibit large spatial variations. When two inhomogeneities are close to each other, intensive interfaces lead to drastic change of heat flux, see Fig. \ref{fig:thermal_spheres} (b) in Section 4 ($x_3 / a_3 \approx 0$). In addition, Figs. \ref{fig:contour_flux} (a-b) show that heat flux decreases with height, including matrix and inhomogeneities regions. It is natural to predict that the averaged heat flux should decrease with height, however, Fig. \ref{fig:AVG_0.8s} (b) and Fig. \ref{fig:AVG_1.5s} (b) provides contradictory conclusions, as the averaged heat flux exhibit fluctuations. The different observation is caused non-representative sampling points. Despite that sampling points are uniformly distributed in one layer (i) when the inhomogeneities are small, only a few points within the inhomogeneity are considered; (ii) when the inhomogeneities are large, the drastic change of heat flux between inhomogeneities are not given full consideration. The above two reasons explains that even adding more layers of inhomogeneities, the fluctuations of averaged heat flux cannot be handled.  Despite inhomogeneities exhibit certain size effects, the heat flux distribution are similar between div = 20 and div = 36, which agrees with Figs. \ref{fig:AVG_0.8s} (b), \ref{fig:AVG_1.5s} (b), and \ref{fig:AVG_3s} (b).

\section{Conclusions}
This paper proposes the dual-reciprocity boundary integral equations (DR-BIEs) for a bi-layered composite system, in which the continuity conditions of temperature and heat flux are taken into account by the bimaterial Green's function. The heat generation/storage rate is considered a nonhomogeneous heat source, whose influences can be handled through two boundary integrals. Eshelby's equivalent inclusion method (EIM) has been extended to simulate inhomogeneities under general transient/time-harmonic heat transfer, and closed-form Eshelby's tensors handle the disturbances caused by eigen-fields. Two eigen-fields, namely eigen-temperature-gradient (ETG) and eigen-heat-source (EHS), are determined by equivalent heat flux and equivalent heat generation conditions, respectively. The combination of DR-BIEs and EIM, namely DR-iBEM, enables the simulation of general heat transfer in bi-layered composites without any interior mesh, including inhomogeneities and bimaterial interface. The DR-BIEs and DR-iBEM algorithms have been validated by convergent FEM results with significant computational advantages in the usage of RAM and efficiency. The DR-iBEM is applied to simulate transient/time-harmonic heat transfer in a functionally graded material (FGM) by switching the bi-layered matrix and inhomogeneities. The gradation effects of inhomogeneities on the FGM are evaluated by the DR-iBEM, which demonstrates its capability in advanced material modeling.

\section*{‌CRediT Author Contribution}
\textbf{Chunlin Wu}: Conceptualization, Methodology, Data Curation, Software, Validation, Writing-Original Draft, Funding Acquisition; \textbf{Liangliang Zhang}: Writing-Review \& Editing; \textbf{Tengxiang Wang}: Writing-Review \& Editing; \textbf{Huiming Yin}: Writing-Review \& Editing, Supervision.

\section*{Acknowledgment}
CW's work is supported by the Chenguang Program No. 23CGA50 of Shanghai Education Development Foundation and Shanghai Municipal Education Commission, National Natural Science Foundation of China Grant No. 12302086; LZ's work is supported by Chinese Universities Scientific Fund 2025TC014; HY's work is sponsored by U.S. Department of Agriculture NIFA \#2021-67021-34201, and NIFA SBIR \#20233353039686, and National Science Foundation (NSF) grants IIP \#1738802 and IIP \#1941244. Those supports are gratefully acknowledged.

\bibliographystyle{unsrt}

\end{document}